\documentclass[12pt]{amsart}
\usepackage{amsmath,amssymb,amscd}

\newtheorem{thm}{Theorem}[section]
\newtheorem{cor}[thm]{Corollary}
\newtheorem{lem}[thm]{Lemma}
\newtheorem{defn}[thm]{Definition}
\newtheorem{prop}[thm]{Proposition}
\newtheorem{examp}[thm]{Example}

\newtheorem{rem}[thm]{Remark}

                \def \x{\mathbf {x}}
              \def \M{\mathcal {M}}
              \def \C{\mathbf {C}}
              \def \Z{\mathbf {Z}}
              \def \R{\mathbf { R}}
              \def \Q{\mathbf { C}}
              \def \U{\mathcal {U}}

              \def \z{\mathbf {z}}
              \def \m{\mathbf { m}}
              \def \n{\mathbf {n}}
              \def \g{\mathbf { g}}
              
              \def \V{\mathcal {V}}

              \def \C{\mathbf {C}}

          \def \H{\mathcal {H}}
 \def \G{{\mathcal G}}

\begin{document}
\noindent
\begin{center}
  {\LARGE Gerbes and twisted orbifold quantum cohomology}\footnote{the second author is partially
supported by the National Science Foundation and Hong Kong RGC
grant}
  \end{center}

  \noindent
  \begin{center}

    {\large Jianzhong Pan, Yongbin Ruan and Xiaoqin Yin}

              \end{center}

\section{Introduction}

    An important part of stringy orbifold theory is the various
    twistings the theory possesses. Unfortunately, it is also the part of
    stringy orbifold theory that we understand the least. For example,
    for the untwisted theory, we have a rather complete
    conjectural
    answer to its structure and relation to crepant resolution. On
    the other hand, both the structure of twisted theory and
    its relation to desingularization are still mysterious at this
    moment. This article fills in one piece of the puzzle.

    Recall that for any almost complex orbifold $X$, we can
    associate a  Chen-Ruan orbifold cohomology ring
    $H^*_{CR}(X, \C)$ \cite{CR1} as the summation of ordinary cohomology
    over all the sectors with appropriate degree shifting. There
    are two important factors of this ring. (1) There is a
    K-theoretic counterpart $K_{orb}(X)$ due to Adem-Ruan \cite{AR}.
    (2) A precise relation between Chen-Ruan orbifold cohomology ring
    and the cohomology ring of its crepant resolution has been
    proposed \cite{R1}. These two key aspects of the theory will always serve as benchmarks to
    our future constructions. Namely, any theory we constructed
    should have two properties: (1) it should be compatible with K-theory; (2) it should describe
    the ring structure of its crepant resolution or more generally desingularization.
    Here, a desingularization is obtained by first deforming the equation of a Gorenstein orbifold
    and then taking a crepant resolution. The miracle is that the right answer for one is often
    automatically the right answer for the other one. This gives us
    two powerful approaches to stringy orbifolds.

    Historically, the earliest twisting is due to Vafa \cite{V}, \cite{VW} in
    the case of a global quotient orbifold $X=Y/G$. Vafa's twisting
    is     a group cohomology class $\alpha\in H^2(G, S^1)$ called {\em discrete torsion}. The notion of
    discrete torsion  was  generalized to arbitrary orbifolds as a class in
    $H^2(\pi^{orb}_1(X), S^1)$ \cite{R}. One can construct a twisted orbifold cohomology using
    discrete torsion \cite{VW}, \cite{R}. Its K-theoretic counterpart was
    constructed by Adem-Ruan \cite{AR}. However, it fails badly
    of describing the cohomology of
    desingularization. To remedy the situation, a more general twisting was proposed
    by the author \cite{R}. This new twisting is a flat line bundle $\mathcal{L}$ over
    the inertial orbifold satisfying a certain compatibility condition
    (see Definition 3.1). $\mathcal{L}$ is called {\em an inner local
    system}. The inner local system works well for the second
    task, i.e, describing the cohomology group of a
    desingularization. But Adem-Ruan's construction of twisted
    orbifold K-theory fails to cover the case of an inner local
    system.

    A more important problem is to twist orbifold quantum
    cohomology, which is unknown even for discrete torsion. Recall
    that the Cohomological Crepant Resolution Conjecture can be
    phrased as follows: the cup product of a crepant resolution $Y$ of
    Gorenstein orbifold $X$ is the Chen-Ruan product of $X$ and
    quantum corrections coming from Gromov-Witten invariants of
    exceptional rational curves. This conjecture was obtained by understanding the
    behavior of quantum cohomology when we deform a crepant resolution to orbifold,
    even though our initial goal is only to understand cohomology. Therefore, to understand even
    the ordinary ring structure of desingularization, one has to
    understand the quantum cohomology and its twisting. This is
    the main goal of this paper.

    Back then, both problems seemed to
    be hopeless. The situation changed when Lupercio-Uribe
    introduced the notion of gerbes to orbifolds \cite{LU1} (see also Tu-Xu-Laurent-Gengoux \cite{TXL}).
    Lupercio-Uribe-Tu-Xu-Laurent-Gengoux constructed twisted K-theory
    using a gerbe on any groupoid, which is much more general than
    an orbifold. Their twisted K-theory generalizes
    Adem-Ruan's twisted K-theory on orbifolds and twisted K-theory
    on smooth manifolds studied by Witten \cite{W}, Bouwknegt-Mathai \cite{BM}
    Freed-Hopkins-Teleman \cite{FHT} and others. The
    beauty of gerbes is that one can easily do differential
    geometry, which is precisely what we were doing for quantum
    cohomology. In this context, Lupercio-Uribe interpreted
    an inner local system as the holonomy line bundle on the inertial
    groupoid of a gerbe. In this article, we would like to go one more
    step further
    to use the gerbe to twist orbifold quantum cohomology.
    During the course of this work, some subtleties arise. In the theory of gerbes, there is a distinction between
    flat gerbes and non-flat gerbes. A flat gerbe has torsion characteristic class and is often referred as
    a torsion gerbe. It has a rather long history in classical geometry under the name of Brauer group. The flat gerbe or
    element of the Brauer group is precisely the data we are able to use to twist orbifold quantum cohomology.
    On a smooth manifold, our twisted orbifold qunatum cohomology did not give any new information (see Corollary \ref{T:smooth}).
    However, an orbifold flat gerbe naturally contains discrete torsion. It gives an abundance of new invariants.
    Our construction does not work for non-flat gerbes. In many ways, non-flat gerbes seem to fall into the
    realm of non-commutative geometry. A further understanding of twisted orbifold quantum cohomology may require a
    full-fledged theory of non-commutative quantum cohomology. The on-going development of geometry with B-field by Hitchin and
    others may provide another approach to this type of question.

      Since a gerbe and its twisted
    K-theory can be constructed over a singular space much more
    general than an orbifold, a natural question is: Can we
    construct a (twisted) orbifold (quantum) cohomology for
    a general groupoid such that (1) it agrees with twisted K-theory
    rationally; (2) it describes the cohomology of its
    desingularization?

    The main results of this paper were announced by the second
    author in 2001 at ICM Satellite conference on Stringy
    Orbifolds
    in Chengdu. For various reason, we were distracted by other
    projects. We apologize for such a long delay. During the
    preparation of this paper, we received an article of Lupercio
    and Uribe where there is some overlap between our section 4
    and their paper \cite{LU}.

    The paper is organized as follows. In section 2, we will
    review the basic definitions of orbifold and groupoid. In
    seciton 3, we will review the definition of gerbes and their
    holonomy. In section 4, we will show how to use the holonomy
    of a gerbe to twist orbifold GW-invariants. Some examples will
    be computed in the last section.

    \section{A review of orbifold Gromov-Witten invariants}
        We will review the construction of ordinary orbifold
        Gromov-Witten invariants due to Chen-Ruan \cite{CR2}. We
        will only sketch the main construction and refer the
        detail to \cite{CR2}. But we take this opportunity to steamline the
    definition.

       From now on, we will use $X_o$ to denote a connected component
    of $X$. We will also assume that all intersections
    $U_{i_1\cdots, i_k}=U_{i_1}\cap \dots \cap U_{i_k}$ are
    connected. Otherwise, we work component by component.

    An orbifold atlas is defined in the same way.
    \vskip 0.1in
    \noindent
    \begin{defn} [Orbifold Atlas] A $n$-dimensional
    smooth orbifold atlas on connected open cover
    $\{U_i\}_{i\in I}$ of $X$ is given by the following data:
    \begin{description}
    \item[(1)]{\it Each $U_i$ is covered by an unformizing system
    $(\tilde{U}_i, G_{\tilde{U}_i}, \pi_i)$ in the following sense.
    $\tilde{U}_i$ is smooth. $G_{\tilde{U}_i}$ is a finite group
    acting smoothly on $\tilde{U}_i$ and $\pi_i:
    \tilde{U}_i\rightarrow U_i$ is invariant under
    $G_{\tilde{U}_i}$ such that it induces a homeomorphism
    $\tilde{U}_i/G_{\tilde{U}_i}\cong U_i$. We call $U_i$ an
    orbifold chart.}

    Choose a component $\pi^{-1}_i(U_{ij})_o$ and let
    $G_{\pi^{-1}_i(U_{ij})_o}\subset G_{\tilde{U}_i}$ be the subgroup fixing $\pi^{-1}_i(U_{ij})_o$.
    Then $(\pi^{-1}_i(U_{ij})_o, G_{\pi^{-1}_i(U_{ij})_o},\pi_i)$
    is an uniformizing system of $U_{ij}$ (called an induced
    uniformizing system). Other induced uniformizing system are of
    the form $(g\pi^{-1}_i(U_{ij})_o,
    gG_{\pi^{-1}_i(U_{ij})_o}g^{-1}, \pi_i)$ for $g\in
    G_{\tilde{U}_i}$. Namely, $G_{\tilde{U}_i}$ acts transitively
    on the collection of induced unformizing systems. In the same
    way, $(\tilde{U}_j, G_{\tilde{U}_j}, \pi_j)$ induces a
    collection of uniformizing systems over $U_{ij}$ acts
    transitively by $G_{\tilde{U}_j}$.

    \item[(2)]{\it $Tran(U_i, U_j)$
    is a collection of isomorphisms from \newline $(\pi^{-1}_i(U_{ij})_o, G_{\pi^{-1}_i(U_{ij})_o},\pi_i)$
    to $(\pi^{-1}_j(U_{ij})_o, G_{\pi^{-1}_j(U_{ij})_o},\pi_j)$.
    Here, the isomorphism is an isomorphism
    $\lambda_{ij}: G_{\pi^{-1}_i(U_{ij})_o}\rightarrow G_{\pi^{-1}_j(U_{ij})_o}$ and an equivariant
    diffeomorphism $\phi_{ij}:\pi^{-1}_i(U_{ij})_o\rightarrow \pi^{-1}_j(U_{ij})_o$.
    Moreover, all such isomorphisms are generated from a fixed one
    by the action of $G_i\times G_j$ in an obvious way. Each
    isomorphism is called a transition map.
    \item[(3)] $Tran(U_i, U_i)$ is generated by the identity.
    $Tran(U_j, U_i)=Tran(U_i, U_j)^{-1}$ in the sense that each
    transition in $Tran(U_j, U_i)$ is the inverse of some
    transition of $Tran(U_i, U_j)$.}

    Over the triple intersection $U_{ijk}$, each of $U_i,U_j, U_k$
    induces a uniformizing system and transitions restrict to them
    as well. Then, we require

    \item[(4)]{\it There is a multiplication such that $(\phi_{jk},
    \lambda_{jk})\circ (\phi_{ij}, \lambda_{ij})$ is the
    restriction of an element of $Tran(U_i, U_k)$.}

    \end{description}
    \end{defn}
    \vskip 0.1in

    Note that we do not require $\tilde{U}_i$ to be connected.

    If $\U'$ is a refinement of $\U$ satisfying (*), then there is
    an induced orbifold atlas over $\U'$ in an obvious fashion. Two orbifold
    atlases are considered to be equivalent if their induced orbifold atlases are
    equivalent over a common refinement in an obvious fashion.  Such an
    equivalence class is called {\em an orbifold structure over
    $X$}. So we may choose $\U$ to be arbitrarily fine.

    Let $x \in X$. By choosing a small neighborhood $V_p\in \U$,
    we may assume that its uniformizing system $\V(V_p)=(U_p, G_p)$
    has the property that $U_p$ is an $n$-ball centered at origin $o$ and
    $\pi_p^{-1}(p)=o$ where $\pi_p: U_p\rightarrow V_p=U_p/G_p$ is the
    projection map. In particular, the origin $o$ is fixed by $G_p$. We
    called $G_p$ {\em the local group at $p$}. If
    $G_p$ acts effectively for every $p$, we call $X$ an {\em effective orbifold}.

    Recall Satake's definition of orbifold map. A map $f:
X\rightarrow Y$ is an orbifold map iff locally $f: U_i\rightarrow
V_i$ can be lifted to an equivariant map $\tilde{f}_i:
\tilde{U}_i\rightarrow \tilde{V}_i$ with a homomorphism
$\lambda_i: G_{\tilde{U_i}}\rightarrow G_{\tilde{V_i}}$. Suppose
that we want to pull back an orbifold vector bundle from $Y$. We
can use a local lifting $\tilde{f}_i$ to construct the local
pull-back. But there is no reason that local pull-backs can be
glued together. In order to glue them together, we have to impose
the condition on transitions. Then, we obtain the notion of {\em
orbifold morphism}. Now it is clear how we should impose our
condition called {\em compatible system}.

    \vskip 0.1in
    \noindent
 \begin{defn}  Fix an underlying map $f:
X\rightarrow Y$. A compatible system consists of an orbifold atlas
$(\tilde{U}_i, G_{\tilde{U}_i}), Tran(U_i, U_j)$ of $X$ and an
orbifold atlas $(\tilde{V}_i, G_{\tilde{V}_i}), Tran(V_i, V_j)$ of
$Y$ with the following additional properties:

\begin{description}
\item[(i)] $f$ maps a member of one atlas to a member of other
atalas, i.e.,$f: U_i\rightarrow V_{\kappa (i)}$. \item[(ii)] The
local map in (i) can be lifted to $\lambda_i:
G_{\tilde{U}_i}\rightarrow G_{ \tilde{V}_{\kappa(i)}}$ and an
equivariant map $\tilde{f}_i: \tilde{U}_i\rightarrow
\tilde{V}_{\kappa(i)}$. \item[(iii)] There is a map $\lambda_{ij}:
Tran(U_i, U_j)\rightarrow Tran(V_{\kappa (i)}, V_{\kappa(j)})$
preserving the
               identity, inverse and multiplication.
\item[(iv)] $\lambda_{ij}(g)\circ \tilde{f}_i=\tilde{f}_j\circ g$.
\end{description}
\end{defn}
\vskip 0.1in

   Suppose that $\{V_{\beta}\}$ is a refinement of $\{V_j\}$.
   Then, $\{f^{-1}(V_{\beta})\}$ is a refinement of $\{U_i\}$. We can
   take a further refinement $\{U_{\alpha}\}$ of
   $\{f^{-1}(V_{\beta})\}$. Then we still have property (i).
   Furthermore, the original compatible system induces
   compatible systems over $\{U_{\alpha}\}$, $\{V_{\beta}\}$.
   We call this a refinement of compatible systems.

   \vskip 0.1in
   \noindent
  \begin{defn} [Isomorphism of compatible systems]

   Two compatible systems given by $(\tilde{f}_i, \lambda_i, \lambda_{ij}),
(\tilde{f}'_i, \lambda'_i, \lambda'_{ij})$ over the same orbifold atlas
   $(\tilde{U}_i, G_{\tilde{U}_i}), (\tilde{V}_j,
   G_{\tilde{V}_j})$ are said to be isomorphic if they differ by
   an automorphism of orbifold structure $(\tilde{V}_j,
   G_{\tilde{V}_j})$. Namely, there is an element $\delta_i\in Tran(V_i, V_i)$
such that
    $$\tilde{f}'=\delta_i\circ \tilde{f}_i,\ \lambda'_i=\delta_i \lambda_i\delta^{-1}_i,
    \ \lambda'_{ij}=\delta_j \lambda_{ij}\delta^{-1}_i.$$
    For two arbitrary compatible systems over isomorphic orbifold atlases, by
taking refinements and components if necessary, we can induce two
compatible systems over the same orbifold atlas. Then the original
ones are isomorphic iff the induced ones are isomorphic in the
above sense. An orbifold morphism is a map with an isomorphism
class of compatible system.
\end{defn}
    \vskip 0.1in
    Chen-Ruan also developed a machinery to classify good maps. The
    key is an invariant they called the characteristic. The case we
    will use is the global quotient orbifold denoted by the stacky
    notation $[X/G]$. The characteristic can be interpreted as follows.

    Suppose that $f: Y\rightarrow [X/G]$ is a good map. We can
    pull back the $G$-bundle $X\rightarrow X/G$ to obtain a $G$-bundle
    $p: E\rightarrow Y$ and a $G$-map $F: E\rightarrow X$. In
    fact, the equivalence class of a good map $f$ is equivalent to
    the pair $(p,F)$ modulo bundle isomorphism. Namely,
    $(p,F)~(p',F')$ iff $p'=ph, F'=Fh$ for a bundle isomorphism
    $h: E'\rightarrow E$.

    Since $G$ is a finite group, $p: E\rightarrow Y$ is an
    orbifold cover. By covering space theory, $E$ is
    determined by the conjugacy class of a homomorphism $\rho:
    \pi^{orb}_1(Y, x_0)\rightarrow G$. We call $\rho$ and its
    conjugacy class Chen-Ruan characteristic.

    Consider the pairs:
$$
\wedge X=\{(p,(g)_{G_p})|p\in X, g\in G_p\},
$$
where $(g)_{G_p}$ is the conjugacy class of $g$ in $G_p$. If there
is no confusion, we will omit the subscript $G_p$ to simplify the
notation. $\wedge X$ has a natural orbifold structure (Proposition
\ref{T:orbi}) and is called {\em the inertia orbifold}. More
generally, we can define {\em the multisector}
$$\widetilde{X}_k=\{(p,(g_1,\cdots, g_k)_{G_p})|p\in X, g_i\in G_p\},
$$
It is clear that $\wedge X=\widetilde{X}_1$. There are two classes
of maps between multisectors.
$$I: \widetilde{X}_k\rightarrow \widetilde{X}_k $$
by
$$I(p,(g_1,\cdots, g_k)_{G_p})=(p,(g^{-1}_1,\cdots,
g^{-1}_k)_{G_p}),$$
    and
    $$e_{i_1, \cdots, i_l}: \widetilde{X}_k \rightarrow
    \widetilde{X}_l $$
    by
    $$e_{i_1, \cdots, i_l}(p,(g_1,\cdots, g_k)_{G_p})=(p,(g_{i_1},
    \dots, g_{i_l})_{G_p}).$$

Suppose that $X$ has an orbifold structure $\U$ with orbifold
atlas $(\tilde{U}_i, G_{\tilde{U}_i}, \pi_i), Tran(U_i, U_j)$.
\vskip 0.1in
\begin{prop} \label{T:orbi}{\it $\widetilde{X}_k$ is naturally an
orbifold, with the orbifold atlas given by
$$ (\bigsqcup_{\g\in G^k_U}\tilde{U}^\g, G_U),
$$
where $\tilde{U}^\g=\tilde{U}^{g_1}\cap \tilde{U}^{g_2}\cap\cdots\cap \tilde{U}^{g_k}$.
 Here
$\g=(g_1,\cdots,g_k)$, $\tilde{U}^g$ stands for the fixed-point
set of $g$ in $\tilde{U}$. When $X$ is almost complex,
$\widetilde{X}_k$ inherits an almost complex structure from $X$,
and when $X$ is closed, $\widetilde{X}_k$ is a finite disjoint
union of closed orbifolds.}
\end{prop}

Next, we would like
   to describe the connected components of $\widetilde{X}_k$. Suppose that
    $p,q$ are in the same orbifold chart $U_i$ uniformized by $(\tilde{U}_i,G_{\tilde{U}_i},\pi_i)$.
 Let $\tilde{p}, \tilde{q}$ be a preimages of $p,q$ respectively. Then $G_p=G_{\tilde{p}}, G_{q}=G_{\tilde{q}}$
and both of them are subgroups of $G_{\tilde{U_i}}$. We say that
$(\g_1)_{G_p}\cong (\g_2)_{G_q}$ if
    $h(q)=p,\g_1=h \g_2 h^{-1}$ for some element $h\in G_{\tilde{U_i}}$. For two arbitrary points $p, q\in X$,
    we say $(\g)_{G_p}\cong (\g')_{G_q}$ if there is a sequence $(p_0, (\g_0)_{G_{p_0}}),
    \cdots, (p_k, (\g_k)_{G_{p_k}})$ such that $(p_0, (\g_0)_{G_{p_0}})=(p, (\g)_{G_p}),
    (p_k, (\g_k)_{G_{p_k}})=(q, (\g')_{G_q})$ and $p_i, p_{i+1}$ are in the same orbifold chart and
    $(\g_i)_{G_{p_i}}\cong (\g_{i+1})_{G_{p_{i+1}}}$. This defines an equivalence relation
    on $\{(\g)_{G_p}\}$. In particular, it is possible that $(\g)_{G_p}\cong (\g')_{G_p}$ while
    $(\g)_{G_p}\neq (\g')_{G_p}$.
  Let $T_k$ be
   the set of equivalence classes.  By  abuse of the notation, we
   often use $(\g)$ to denote the equivalence class to which $(\g)_{G_q}$
   belongs.
   It is clear that $\widetilde{X}_k$ decomposes as a disjoint
   union of connected components
   $$
   \widetilde{X}_k=\bigsqcup_{(\g)\in T_k} X_{(\g)},
   $$
   where
   $$
   X_{(\g)}=\{(p,(\g')_{G_p})|\g'\in G^k_p, (\g')_{G_p}\in (\g)\}.
   $$
    Let $T^o_k\subset T^k$ be such that $(g_1, \cdots, g_k)\in T^o_k$ has the property $g_1\cdots g_k=1$.

    $$\overline{\M}_k(X)=\bigsqcup_{(\g)\in T^o_k} X_{(\g)}.$$
   \vskip 0.1in

   \begin{defn}{\it $X_{(g)}$ for $g\neq 1$ is called
   a twisted sector.  $X_{(\g)}$ is called a $k$-multi-sector or a $k$-sector.
   Furthermore, we call $X_{(1)}=X$ the nontwisted sector.}
   \end{defn}
   \vskip 0.1in

\vskip 0.1in \begin{examp}{\it Suppose that $X=Y/G$ is a global
quotient. By the definition, $\wedge X=\bigsqcup_{g\in G} Y^g/G$
where $Y^g$ is the fixed-point set of elements $g\in G$.
Equivalently, $\wedge X=\bigsqcup_{(g)} Y^g/C(g)$.}
\end{examp}
 \vskip 0.1in
\begin{prop}{\it Both the evaluation maps $e_{i_1, \cdots,
i_l}$ and $I$ are orbifold morphisms.}
\end{prop}
\vskip 0.1in

    Next, we extend the notion of orbifold morphism to the case in which the
    domain is a nodal orbifold Riemann surface.

    Recall that
a  nodal curve with $k$ marked points is a pair $(\Sigma,\z)$
consisting of a connected topological space
$\Sigma=\bigcup\pi_{\Sigma_\nu}(\Sigma_\nu)$, where $\Sigma_\nu$
is a smooth complex curve  and $\pi_\nu:\Sigma_\nu\rightarrow
\Sigma$ is a continuous map, and $\z=(z_1,\cdots,z_k)$ consists of
$k$ distinct  points in $\Sigma$ with the following properties:
\begin{itemize}
\item For each $z\in\Sigma_\nu$, there is a neighborhood of it such that
the restriction of $\pi_\nu:\Sigma_\nu\rightarrow \Sigma$ to this
neighborhood is a homeomorphism to its image.
\item For each $z\in\Sigma$, we have $\sum_\nu \#\pi_\nu^{-1}(z)\leq 2$.
\item $\sum_\nu \#\pi_\nu^{-1}(z_i)=1$ for each $z_i\in\z$.
\item The number of complex curves $\Sigma_\nu$ is finite.
\item The set of nodal points $\{z|\sum_\nu \#\pi_\nu^{-1}(z)=2\}$ is  finite.
\end{itemize}

A point $z\in\Sigma_\nu$ is called {\em singular} (or a
\emph{node}) if $\sum_\omega\#\pi_\omega^{-1}(\pi_\nu(z))=2$. A
point $z\in\Sigma_\nu$ is said to be a {\em marked point} if
$\pi_\nu(z)=z_i\in\z$. Each $\Sigma_\nu$ is called a {\em
component} of $\Sigma$. Let $k_\nu$ be the number of points on
$\Sigma_\nu$ which are either singular or marked, and $g_\nu$ be
the genus of $\Sigma_\nu$; a nodal curve $(\Sigma,\z)$ is called
{\em stable} if $k_\nu+2g_\nu\geq 3$ holds for each component
$\Sigma_\nu$ of $\Sigma$.

    \vskip 0.1in
    \noindent
\begin{defn}{\it \emph{A nodal orbicurve} is a nodal marked  curve
$(\Sigma,\z)$ with an orbifold structure as follows:
\begin{itemize}
\item The set $\z_\nu$ of orbifold points of each
component $\Sigma_\nu$ is contained in the set of marked points
and nodal points $\z$.
\item A neighborhood of a marked point
is uniformized by a branched covering map $z\rightarrow z^{m_i}$
with $m_i\geq 1$.
\item A neighborhood of a nodal point
(viewed as a neighborhood of the origin of $\{x y=0\}\subset
\C^2$) is uniformized by a branched covering map $(x,y)\rightarrow
(x^{n_j}, y^{n_j})$, with $n_j\geq 1$, and with group action $e^{2
\pi i /n_j}(x,y)=(e^{2 \pi i /n_j}x, e^{- 2\pi i/n_j}y)$.
\end{itemize}
Here $m_i$ and $n_j$ are allowed to be equal to one, i.e., the
corresponding orbifold structure is trivial there. We denote the
corresponding nodal orbicurve by $(\Sigma,\z,\m,\n)$ where
$\m=(m_1,\cdots,m_k)$ and $\n=(n_j)$.}
\end{defn}
    \vskip 0.1in
    Once we have the definition of nodal orbicurve, we can extend
    the definition of compatible system and orbifold morphism word by word to the case where the
    domain is a nodal orbicurve.

    First, recall that for every point $p \in \Sigma$, an orbifold morphism
 $f: \Sigma \rightarrow X$ induces a homomorphism $G_p
\rightarrow G_{f(p)}$.

    \vskip 0.1in
    \noindent
  \begin{defn}{\it Let $(X,J)$ be an almost complex orbifold. An \emph{
orbifold stable map into $(X,J)$} is a quadruple
$(f,(\Sigma,\z,\m,\n), \varphi_s), \xi)$ described as follows:

\begin{enumerate}
\item $f$ is a continuous map from the nodal orbicurve
$(\Sigma,\z,\m,\n)$ into $X$ such that each $f_\nu=f\circ\pi_\nu$
is a pseudo-holomorphic map from $\Sigma_\nu$ into $X$. \item
$\xi$ is an isomorphism class of compatible structures. \item Let
$k_\nu$ be the order of the set $\z_\nu$, namely the number of
points on $\Sigma_\nu$ which are special (i.e. nodal or marked );
if $f_\nu$ is a constant map, then $2 g_\nu - 2  + k_\nu
> 0$.
\item At any marked or nodal point $p$ the induced homomorphism on
the local group $\lambda_{p}: G_{p}\rightarrow G_{f(p)}$ is
injective.
\end{enumerate}}
\end{defn}

Finally we observe that each $C^\infty$ orbifold morphism
 from an orbifold nodal Riemann surface with $k$ marked
points into an orbifold $X$ determines a point in the product of
inertia orbifolds $(\wedge X)^k$ as follows: let the underlying
continuous map be $f$ and for each marked point $z_i$,
$i=1,\cdots,k$, let $x_i$ be the positive generator of the cyclic
local group at $z_i$, and $\lambda_{z_i}$ be the homomorphism
determined by the given compatible system; then the determined
point in $(\wedge X)^k$ is
$$
((f(z_1),(\lambda_{z_1}(x_1))_{G_{f(z_1)}}),\cdots,
(f(z_k),(\lambda_{z_k}(x_k))_{G_{f(z_k)}})).
$$
Let $\x=(X_{(g_1)},\cdots,X_{(g_k)})$ be a connected component in
$(\wedge X)^k$. We say that a good map with a compatible system is
of {\it type} $\x$ if the above point  it determines in $(\wedge
X)^k$ lies in the component $\x$.

\vskip 0.1in

\begin{rem}{\it If $f: \Sigma\rightarrow X$ is a
pseudo-holomorphic map whose image intersects the singular locus
of $X$ at only finitely many points, then there is a unique choice
of orbifold structure on $\Sigma$ together with a unique
$(\tilde{f},\xi)$, where $\tilde{f}$ is a good map with an
isomorphism class of compatible systems $\xi$ whose underlying
continuous map is $f$. If the image of $f$ lies completely inside
the singular locus, there could be different choices, and they are
regarded as different points in the moduli space. }
\end{rem}

\vspace{2mm}

\begin{defn}\hspace{2mm}{\it
\begin{enumerate}
\item An orbifold $X$ is symplectic if there is a closed 2-form $\omega$ on
$X$ whose local liftings are non-degenerate.
\item A projective orbifold is a complex orbifold which is a projective
variety as an analytic space.
\end{enumerate}
}
\end{defn}
\vspace{2mm}

\begin{prop}\hspace{2mm}{\it Suppose that $X$ is a symplectic or
projective orbifold. The moduli space of orbifold stable maps
$\overline{\M}_{g,k}(X,J,A, \x)$ is a compact metrizable space
under a natural topology, whose ``virtual dimension'' is $2d$,
where $$ d=c_1(TX)\cdot A+(\dim_{\C}X-3)(1-g)+k-\iota(\x).
$$ Here $\iota(\x):=\sum_{i=1}^k \iota_{(g_i)}$ for
$\x=(X_{(g_1)},\cdots,X_{(g_k)})$. }
\end{prop}

For any component $\x=(X_{(g_1)},\cdots,X_{(g_k)})$, there are $k$
evaluation maps  $$ e_i:\overline{\M}_{g,k}(X,J,A,\x)\rightarrow
X_{(g_i)}, \hspace{4mm} i=1,\cdots, k. $$
    $e_i$ has a natural compatible system to make it a good map. For any set of
cohomology classes $\alpha_i\in H^*(X_{(g_i)};\Q)\subset
H^*_{orb}(X;\Q)$, $i=1,\cdots,k$, the orbifold Gromov-Witten
invariant is defined as  $$
\Psi^{X,J}_{(g,k,A,\x)}(\alpha^{l_1}_1, \cdots, \alpha^{l_k}_k)=
\prod_{i=1}^k c_1(L_i)^{l_i}e^*_i \alpha_i
[\overline{\M}_{g,k}(X,J,A,\x)]^{vir},$$ where $L_i$ is the line
bundle generated by the cotangent space of the $i$-th marked
point. The virtual fundamental cycle
$[\overline{\M}_{g,k}(X,J,A,\x)]^{vir}$ is defined as the
fundamental cycle of a certain orbifold $S$ \cite{CR2}.

\vspace{2mm}

  The inertial orbifold admits another interpretation as the space
  of   constant loops. Then it is naturally a subset of the free loop
  space. We shall sketch this construction due to Chen \cite{C}
  (see \cite {LU2} for a groupoid description).

  Let $\Omega X$ be the space of orbifold morphisms from $S^1$ with trivial
  orbifold structure to $X$. $\Omega X$ is the analog of the free
  loop space of a smooth manifold. We need the following important

  \vskip 0.1in
  \noindent
 \begin{lem}[Lemma 3.15 in \cite{C}]{\it Let $X=Y/G$ be a
  global quotient. Then, $\Omega X=P(Y, G)/G$, where
  $$P(Y, G)=\{(\gamma, g); \gamma: [0,1]\rightarrow Y, g\in G,
  \gamma(1)=g\gamma(0)\}.$$
  Here, $G$ acts on $P(Y, G)$ by $h(\gamma,g)=(h\circ \gamma,
  h^{-1}gh)$.}
  \end{lem}
  \vskip 0.1in
  We call $\tau$ a constant loop if the underlying map is constant.
  Suppose that the image is $p\in X$. Let $U_p/G_p$ be the
  orbifold chart at $p$. By the Lemma, $\Omega U_p/G_p=P(U_p, G_p)/G$.
  In particular, $\tau$ is an equivalence class of a pair $(\gamma, g)$
  where im $(\gamma)=p$. Under the action of $G$, we naturally
  identify it as $(p, (g)_{G_p})$. Therefore, the space of
  constant loop is precisely the inertia orbifold $\widetilde{X}$.

  Suppose that $f: \Sigma \rightarrow X$ is an orbifold stable
  map. We take a real blow-up of $\Sigma$ at all the marked points to obtain a Riemann
  surface with boundary $\Sigma^{\dagger}$. $\Sigma^{\dagger}$ can be understood as follows.
  We remove the marked point $x_i$. A neighborhood of a puncture point $x_i$ is
  biholomorphic to $S^1\times [0, \infty)$. Hence, we can view
  $\Sigma^{\dagger}$ as a manifold with cylindrical end and $x_i$ is replaced by
  a circle $S_{\infty}$ attached at $\infty$. Another way
  to interpret the evaluation map is that
  $$e_i(f)=f(S_{\infty}).$$
    This description is important later in our construction.

  \section{Gerbes and their holonomy}
    After reviewing the construction of orbifold quantum cohomology in the last section,
we are ready to touch upon the main topic of this
article-twisting. The earlist twisting from physics is discrete
torison by Vafa \cite{V}. However, discrete torison is too
restrictive to describe interesting examples. Therefore, a more
general twisting is needed. For this purpose, the second author
introduced the notion of inner local system. Roughly speaking, an
inner local system is a flat orbifold line bundle over the
inertial orbifold $\widetilde{X}_1$ satisfying certain
compatibility conditions. Later, Lupercio-Uribe introduced the
concept of gerbe to orbifolds. The holonomy line bundle of a gerbe
with connection is naturally an inner local system. However, not
all inner local systems are induced in this way \cite{AP}. In this
section, we will study the relation between a gerbe and its
holonomy in detail.
    \subsection{Inner local system}
    Recall that for $(g_1, \cdots, g_k)\in T_k$, there are $k+1$
evaluation maps
$$e_i: X_{(g_1, \cdots, g_k)}\rightarrow X_{(g_i)},\ i\leq k,$$
and
$$ e_{k+1}: X_{(g_1, \cdots, g_k)}\rightarrow X_{(g_1\cdots g_k)}.$$

    Now we introduce the notion of inner local system for an orbifold.
\vskip 0.1in

\begin{defn} Suppose that $X$ is an orbifold (almost complex or not).
An {\em inner  local system} $\mathcal{L}=\{L_{(g)}\}_{g\in T_1}$
is an assignment of  a flat complex orbifold line bundle
$$
L_{(g)}\rightarrow X_{(g)}
$$
to each sector $X_{(g)}$ satisfying the following compatibility
conditions (1-4).
\begin{description}
\item[(1)] $L_{(1)}$ is a trivial orbifold line bundle with a
fixed trivialization. \item[(2)] There is a nondegenerate pairing
$L_{(g)}\otimes I^*L_{(g^{-1})}\rightarrow \C=L_{(1)}.$ \item[(3)]
There is a multiplication
$$e_1^*L_{(g_1)}\otimes e_2^*L_{(g_2)}\stackrel{\theta}{\rightarrow} e^*_3 L_{(g_1g_2)}$$
over $X_{(g_1, g_2)}$ for $(g_1, g_2)\in T_2$. \item[(4)] $\theta$
is associative in the following sense. For $(g_1,g_2, g_3)\in
T_3$, the evaluation maps $e_i: X_{(g_1, g_2, g_3)}\rightarrow
X_{(g_i)}$  factor through
$$P=(P_1,P_2): X_{(g_1, g_2, g_3)}\rightarrow X_{(g_1, g_2)}\times X_{(g_1g_2, g_3)}.$$
Let $e_{12}:  X_{(g_1, g_2,g_3)}\rightarrow X_{(g_1g_2)}$. We
first use $P_1$ to define
$$\theta: e_1^*L_{(g_1)}\otimes e_2^*L_{(g_2)}\rightarrow e^*_{12}L_{(g_1g_2)}.$$
Then, we can use $P_2$ to define a product
$$\theta: e^*_{12}L_{(g_1g_2)}\otimes e^*_3L_{(g_3)}\rightarrow e^*_4 L_{(g_1g_2g_3)}.$$
Taking the composition, we define
$$\theta(\theta(e^*_1 L_{(g_1)},e^*_2L_{(g_2)}), e^*_3L_{(g_3)}):
e^*_1L_{(g_1)}\otimes e^*_2L_{(g_2)}\otimes
e^*_3L_{(g_3)}\rightarrow e^*_4L_{(g_4)}.$$ On the other hand, the
evaluation maps $e_i$  also factor through
$$P': X_{(g_1,g_2,g_3)}\rightarrow X_{(g_1, g_2g_3)}\times X_{(g_2, g_3)}.$$
In the same way, we can define another triple product
$$\theta(e^*_1L_{(g_1)},\theta(e^*_2L_{(g_2)}, e^*_3L_{(g_3)})):
e^*_1L_{(g_1)}\otimes e^*_2L_{(g_2)}\otimes
e^*_3L_{(g_3)}\rightarrow e^*_4L_{(g_4)}.$$ Then, we require the
associativity
$$\theta(\theta(e^*_1 L_{(g_1)},e^*_2L_{(g_2)}), e^*_3L_{(g_3)})=\theta(e^*_1L_{(g_1)},
\theta(e^*_2L_{(g_2)}, e^*_3L_{(g_3)})).$$

\end{description}
\end{defn} If $X$ is a complex orbifold, we assume that $L_{(g)}$ is
holomorphic.

        \vskip 0.1in
    \noindent
   \begin{defn}{\it Given an inner local system $\mathcal{L}$, we define the twisted orbifold cohomology
    $$H^*_{CR}(X, \mathcal{L})=\oplus_{(g)}
    H^{*-2\iota_{(g)}}(X_{(g)}, L_{(g)}).$$}
    \end{defn}
    \vskip 0.1in
    \noindent
   \begin{defn}{\it Suppose that $X$ is a closed complex
    orbifold and $\mathcal{L}$ is an inner local system. We define Dolbeault cohomology groups
    $$H^{p,q}_{CR}(X, \mathcal{L})=\oplus_{(g)}H^{p-\iota_{(g)},
    q-\iota_{(g)}}(X_{(g)}; L_{(g)}).$$}
    \end{defn}
    \vskip 0.1in
    \noindent
  \begin{prop}{\it If $X$ is a K\"{a}hler orbifold, we have the
    Hodge decomposition
    $$H^k_{CR}(X, \mathcal{L})=
\oplus_{k=p+q} H^{p,q}_{CR}(X, \mathcal{L}).$$}
\end{prop}
    \vskip 0.1in
    \noindent
    {\bf Proof: } Note that each sector $X_{(g)}$ is a K\"ahler orbifold.
    The proposition follows by applying the ordinary Hodge theorem with twisted
    coefficients to each sector $X_{(g)}$. $\Box$

    \subsection{Basics on gerbes and connections}
    The original motivation for the  introduction of gerbes to orbifolds by Lupercio-Uribe is to
understand inner local systems conceputally. Let's start from the
definition of a gerbe on a smooth manifold. We follow closely the
exposition of \cite{H}.

    Let's  suppose $X$ is a smooth manifold and $\{U_\alpha
\}_\alpha$ an open cover. Recall the definition of line bundle. It
can be described by transition functions
    $$g_{\alpha\beta}: U_{\alpha\beta}=U_{\alpha}\cap U_{\beta}\rightarrow S^1$$
satisfying the conditions
    $$g_{\alpha\alpha}=1, g_{\beta\alpha}=g_{\alpha\beta}^{-1}, (\delta g)_{\alpha\beta\gamma}
    =g_{\alpha}g_{\beta}g_{\gamma}=1.$$
In terms of  cohomological language, $g_{\alpha\beta}$ is a
$\check{C}ech$ 1-cocycle of the sheaf of $S^1$-valued functions
$C^{\infty}(S^1)$. Two sets of transition functions induce
isomorphic line bundles iff they induce the same class in $H^1(X,
C^{\infty}(S^1))$.

   A gerbe is a generalization of a line bundle. It is defined as a $\check{C}ech$ 2-cocycle of  sheaf of $S^1$-valued
function $C^{\infty}(S^1)$ over some open cover $\U$. Two gerbes
are {\em equivalent} if they induced the same cocycle over a
common refinement. They are isomorphic if they induced the same
cohomology class in $H^2(X, C^{\infty}(S^1))$. Let
$\U=\{U_{\alpha}\}$ be an open cover. In terms of local data, they
are functions
$$g_{\alpha \beta \gamma} : U_\alpha \cap U_\beta \cap U_\gamma \rightarrow S^1$$
defined on the threefold intersections satisfying
$$g_{\alpha \beta \gamma} = g_{\alpha \gamma \beta}^{-1} = g_{\beta \alpha \gamma}^{-1} = g_{\gamma \beta \alpha}^{-1}$$
and the cocycle condition
$$(\delta g)_{\alpha \beta \gamma \eta} = g_{\beta \gamma \eta} g_{\alpha \gamma \eta}^{-1}
g_{\alpha \beta \eta} g_{\alpha \beta \gamma}^{-1} =1 $$ on the
four-fold intersections $U_\alpha \cap U_\beta \cap U_\gamma \cap
U_\eta$.
 It
also defines a class in $H^3(X; \Z)$; Consider the long exact
sequence of cohomology
$$\cdots \rightarrow H^i(X,C^\infty(\R)) \rightarrow H^i(X,C^\infty(S^1)) \stackrel{\tau_i}{\rightarrow} H^{i+1}(X, \Z) \rightarrow \cdots$$
 derived from the exact sequence of sheaves
$$ 0 \rightarrow \Z \rightarrow C^\infty(\R) \rightarrow C^\infty(S^1) \rightarrow 1.$$
Recall that $\tau_1([g_{\alpha\beta}])\in H^2(X, \Z)$ is the first
Chern class of the corresponding line bundle. In the same way, the
characteristic class of a gerbe is
$\tau_2([g_{\alpha\beta\gamma}])$. It is well-known that  $
C^{\infty}(\R)$ is a fine sheaf;  we get $H^2(X, C^\infty(S^1))
\cong H^3(X,\Z)$. We might say that a gerbe is determined
topologically by its characteristic class. Furthermore, we can
tensor them using the product of cocycles.

 We call a gerbe $g=\{g_{_{\alpha\beta\gamma}}\}$ {\em a trivial gerbe} if $g=\delta f$ is a coboundary for some 1-cochain $f$. $f$
is called a trivialization of $g$. In terms of local data, $f$ is defined by functions
$$f_{\alpha \beta} =f_{\beta \alpha} : U_\alpha \cap U_\beta \rightarrow S^1$$
on the twofold intersections such that
 $$g_{\alpha \beta \gamma} = f_{\alpha \beta} f_{\beta \gamma} f_{\gamma \alpha}$$
Hence, $g$ is represented as a coboundary $\delta f = g$.

Suppose that $f_1, f_2$ are two different trivializations of $g$.
Then $\delta(f_1 f^{-1}_2)=1$. Hence $h=f_1f^{-1}_2$ is a
1-cocycle and hence defines a line bundle.

 A  connection will consist of a pair $(A_{\alpha\beta}, F_{\alpha})$ where
$A_{\alpha\beta}$ are a 1-forms over the double intersections
$A_{\alpha \beta}$, such that
$$iA_{\alpha \beta} + i A_{\beta \gamma} + i A_{\gamma \alpha}
= g_{\alpha \beta \gamma}^{-1} d g_{\alpha \beta \gamma}$$ and the
2-forms $F_\alpha$ are defined over $U_\alpha$ such that $F_\alpha
- F_\beta = dA_{\alpha \beta}$. Note that we define a global
3-form $G$ such that $G|_{U_\alpha} = F_{\alpha}$. This 3-form $G$
is called the curvature of the gerbe connection.

When the curvature $G$ vanishes we say that the connection on the
gerbe is \emph{flat}. Therefore, $dF_{\alpha}=0$. Since
$U_{\alpha}$ is contractible, we can find $B_{\alpha}$ such that
$F_{\alpha}=d B_{\alpha}$. Then, on $U_{\alpha}\cap U_{\beta}$,
    $$F_{\beta}- F_{\alpha}= d A_{\alpha \beta}=d (B_{\beta}-
    B_{\alpha}).$$
    This implies that
    $$A_{\alpha \beta}-B_{\beta}+B_{\alpha}=d f_{\alpha \beta}.$$
    From the definition of connection
    $$iA_{\alpha \beta}+i A_{\beta \gamma} + i A_{\gamma \alpha}=
    g^{-1}_{\alpha \beta \gamma} d g_{\alpha \beta \gamma}.$$
    Hence,
    $$d(if_{\alpha \beta}+i f_{\beta \gamma}+ if_{\gamma\alpha}-\log
    g_{\alpha\beta\gamma})=0.$$
    Let
    $$c_{\alpha \beta\gamma}=e^{if_{\alpha \beta}}e^{i f_{\beta \gamma}}e^{if_{\gamma\alpha}}
    g_{\alpha\beta\gamma}.$$
    $c_{\alpha\beta\gamma}$ is constant. It is clear that $c_{\alpha\beta\gamma}$ is a
    2-cocycle differing from $g_{\alpha\beta\gamma}$ by a coboundary $e^{if_{\alpha \beta}}e^{i f_{\beta \gamma}}e^{if_{\gamma\alpha}}$. Since it is constant,  $c_{\alpha\beta\gamma}$ represents a
    $\check{C}ech$ class in $H^2(X, S^1)$ which we call the
    {\em holonomy} of the connection.

    Next, we check that $\{c_{\alpha\beta\gamma}\}$ is independent of the choice of $B_{\alpha}, f_{\alpha\beta}$ as
    a $\check{C}ech$ cohomology class. This is the analogue of the fact that one can use a flat connection on a line bundle
    to change the transition function to be constant. If we have different $B'_{\alpha}$, then
    $d(B_{\alpha}-B'_{\alpha})=0$ and hence we can write $B_{\alpha}-B'_{\alpha}=df_{\alpha}$.
    Let
    $$f'_{\alpha\beta}=f_{\alpha\beta}-f_{\alpha}+f_{\beta}.$$
    $$df'_{\alpha\beta}=A_{\alpha\beta}-B_{\beta}+B_{\alpha}-B_{\alpha}+B'_{\alpha}+B_{\beta}-B'_{\beta}
    =A_{\alpha\beta}-B'_{\beta}+B'_{\alpha}.$$
    Then
    $$c'_{\alpha\beta\gamma}=e^{if'_{\alpha \beta}}e^{i f'_{\beta \gamma}}e^{if'_{\gamma\alpha}}
    g_{\alpha\beta\gamma}=c_{\alpha\beta\gamma}.$$
    If we have a different choice
    $$A_{\alpha \beta}-B_{\beta}+B_{\alpha}=d \tilde{f}_{\alpha \beta},$$
    $$\tilde{f}_{\alpha\beta}=f_{\alpha\beta}+\lambda_{\alpha\beta}$$
    where $\lambda_{\alpha\beta}$ is a constant function.
    Then,
    $$\tilde{c}_{\alpha\beta\gamma}=e^{i\tilde{f}_{\alpha \beta}}e^{i \tilde{f}_{\beta \gamma}}
    e^{i\tilde{f}_{\gamma\alpha}}    g_{\alpha\beta\gamma}=c_{\alpha\beta\gamma}e^{i\lambda_{\alpha\beta}}
    e^{i\lambda_{\beta\gamma}}e^{i\lambda_{\gamma\alpha}}.$$
    Namely, it differs by a coboundary in the constant sheaf $S^1$.

    For a line bundle, when the holonomy is trivial we get a
    covariant constant trivialization of the bundle.
    If the
    holonomy of a gerbe is trivial, then $c_{\alpha\beta\gamma}$
    is a coboundary, so that there are constants $k_{\alpha
    \beta}\in  S^1$ such that
    $$ c_{\alpha\beta\gamma}=k_{\alpha
    \beta}k_{\beta\gamma}k_{\gamma\alpha}.$$
    Let
    $$h_{\alpha\beta}=k_{\alpha\beta}e^{-if_{\alpha\beta}}.$$
    Then,
    $$h_{\alpha\beta}h_{\beta\gamma}h_{\gamma\alpha}=g_{\alpha\beta\gamma}$$
    and so we have a trivialization of the gerbe, which we call a
    flat trivialization.

    Suppose that the line bundle is given by a $\check{C}ech$ cocycle
    $g_{\alpha\beta}$. Recall that a connection is a 1-form $A_{\alpha}$ on
    $U_{\alpha}$ such that
    $$iA_{\beta}-iA_{\alpha}=g^{-1}_{\alpha\beta}d
    g_{\alpha\beta}.$$
    A section is $f_{\alpha}: U_{\alpha}\rightarrow S^1$ such that
    $f_{\alpha}=g_{\alpha\beta} f_{\beta}$. It is covariant
    constant iff it satisfies the equation
    $df_{\alpha}=iA_{\alpha}f_{\alpha}$. If we write $f_{\alpha}=e^{i
    p_{\alpha}}$, then $dp_{\alpha}=A_{\alpha}$. Therefore, a
    necessary condition is $F_{\alpha}=dA_{\alpha}=0$, i.e., the connection is
    flat. In the case of a gerbe, $dA_{\alpha\beta}\neq 0$ in
    general. We have to allow the freedom to choose $B_{\alpha}$
    such that $d(A_{\alpha\beta}-B_{\beta}+B_{\alpha})=0$. Hence, the
    trivialization $h_{\alpha\beta}$ satisfies a modified equation
    $$dh_{\alpha\beta}=i(A_{\alpha\beta}-B_{\beta}+B_{\alpha})h_{\alpha\beta}.$$

Suppose that we have a second flat
    trivialization $h'_{\alpha\beta}$; then
    $g_{\alpha\beta}=h'_{\alpha\beta}/h_{\alpha\beta}$ defines a
    line bundle $L$. Moreover,
    $$iB_{\beta}-i B_{\alpha}-i A_{\alpha\beta}=d \log
    h_{\alpha\beta},$$
    $$i B'_{\beta}- i B'_{\alpha} -i A_{\alpha\beta}=d\log
    h'_{\alpha\beta}.$$
    Hence,
    $$i(B'-B)_{\beta}-i (B'-B)_{\alpha}=d \log g_{\alpha\beta}$$
    and $A_{\alpha}=(B'-B)_{\alpha}$ defines a connection on $L$.
    By the definition of $B_{\alpha}$ and $B'_{\alpha}$,
    $$F_{\alpha}=d B_{\alpha}=d B'_{\alpha}.$$
    Hence the curvature $dA_{\alpha}=0$. Thus, the difference
    of two flat trivializations of a gerbe is a flat line bundle.
    One can show that the converse is also true.

    \subsection{String connection}

    Recall that a connection on a line bundle induces a holonomy map $Hol: \Omega X\rightarrow S^1$.
The holonomy of a connection on a gerbe has similar property. One
way to understand it is via its analogy to topological quantum
field theory. Recall that  topological quantum field theory can be
described as follows. For any oriented $d$-dimensional manifold
$D$, we associate a Hilbert space $\H_D$. For any cobordism $W$
such that $\partial W=D_1\cup -D_2$, we associate a homomorphism
$\theta_W: \H_{D_1}\rightarrow \H_{D_2}$. $\theta_W$ satisfies the
gluing axiom. Suppose that $\partial W_{12}=D_1\cup -D_2,
\partial W_{23}=D_2\cup -D_3$. We can glue $W_{12}, W_{23}$ along
$D_2$ to obtain $W_{13}$. Then the gluing axiom is
$\theta_{13}=\theta_{23}\circ \theta_{12}$. The analogy for a
gerbe is called {\em a string connection}. It contains the
following ingredients:

   (i) Let $l: S^1\rightarrow X$ be a smooth map. Since $S^1$ is one-dimensional,
    the pull-back of a gerbe with connection to the circle is flat and
    has trivial holonomy. Thus we have flat trivializations. For each $l$, we associate
the moduli space of flat trivializations $\mathcal{L}_l$.
$\mathcal{L}_l$ is analogous to $\H_D$. Recall that we identify
    flat trivializations if they differ by a flat line bundle with trivial holonomy.
        Then, for each loop we have  a space which is acted on freely and transitively by
    the moduli space of flat line bundles $H^1(S^1, S^1)\cong S^1$. Hence, $\mathcal{L}_l$ is isomorphic to
    $S^1$. In other words we
    have a principal $S^1$ bundle $\mathcal{L}$ over the free loop space $\Omega X$.

    We will pay special attention to the space of constant loops.
    Since $X$ is embedded in $\Omega X$ as the space of constant loops,
    it is interesting to compute the restriction of $\mathcal{L}$ over $X$.
    Suppose that $f: S^1\rightarrow X$ is a constant map. Then
    $f$ is the composition of $p:S^1\rightarrow pt$ and $i_f:
    pt\rightarrow X$. The pull-back gerbe
    $i_f^*g_{\alpha\beta\gamma}$ is obviously trivial. Furthermore,
    $i_f^*F_{\alpha}=0, i_f^*A_{\alpha\beta}=0$. Any trivialization of
    $i_f^*g_{\alpha\beta\gamma}$ is a flat trivialization. A key
    observation is that the flat line bundle over a point is
    trivial as well. Therefore, the pull-back gerbe with
    connection by $i_f$ fixes a unique flat trivialization. Its
    pull-back by the projection map $p:S^1\rightarrow pt$ defines a
    canonical element $s_f\in \mathcal{L}_f$.
    and hence a canonical section of $\mathcal{L}|_X$. Hence, $\mathcal{L}|_X$
    is trivial with a canonical trivialization. Therefore, we
    obtain
    \vskip 0.1in
    \noindent
    \begin{lem}{\it $\mathcal{L}|_X$ is independent of the connection of
    the gerbe. Furthermore, it is trivial with a canonical
    trivialization.}
    \end{lem}
    \vskip 0.1in

    (ii) Suppose that $f: \Sigma \rightarrow X$, where $\Sigma$ is a closed Riemann surface.
    Then the pull-back connection of the gerbe $(f^*F_{\alpha},  f^*A_{\alpha,\beta})$ is flat.
    Its holonomy $Hol_f=\{c_{\alpha\beta\gamma}\}$ is a cohomology class in $H^2(\Sigma, S^1)$.
    Since $H^2(\Sigma, S^1)=Hom(H_2(\Sigma, \Z), S^1)$, its evaluation on the fundamental class
    of $\Sigma$ naturally identifies it as a complex number.

    A more interesting case is the case of a Riemann surface with boundaries. Suppose that $f: \Sigma
    \rightarrow X$, where $\Sigma$ is a Riemann surface with
    boundary $t_i$ with a fixed orientation-preserving parameterization
    $\delta_i: S^1\rightarrow t_i$. Let $l_i=f\circ \delta_i$. Since
    $H^2(\Sigma, S^1)=H^3(\Sigma, \Z)=0$, the pull-back gerbe is
    trivial and its holonomy is trivial as well. A flat
    trivialization restricts to a flat trivialization on each
    $t_i$. Namely, it induces an element $\sigma$ in $\prod_i \mathcal{L}_{l_i}$. A
    different flat trivialization of $\Sigma$ differs by a flat
    line bundle of $\tau$ of $\Sigma$. It restricts to a flat line
    bundle $\tau_i$ over each boundary circle viewed as a standard
    $S^1$ via $\delta_i$. Recall that
    $$\pi_1(\Sigma)=\{\lambda_1, \cdots, \lambda_{2g}, l_1,
    \cdots,
    l_k| \prod_i [\lambda_{2i-1}, \lambda_{2i}]
    l_1\cdots l_k=1\}.$$
    Hence, $\tau_1\cdots\tau_k=1$.
    Hence, different flat trivializations induce  elements
    differing by a multiplication of $(\tau_1, \cdots, \tau_k)$ of
    $\tau_i\in S^1$ with $\tau_1\cdots\tau_k=1$.
    Suppose that $L_1, \cdots, L_k$ are k-circle bundles.
    $L_1\otimes \cdots \otimes L_k$ can be constructed as follows.
    Let $H$ be the $(S^1)^k$-bundle with fiber $\prod_i (L_i)_x$.
    Then $\otimes_i L_i=H\times S^1/(S^1)^k$ via the product
    homomorphism $(S^1)^k\rightarrow S^1$. From the previous
    construction,
    \vskip 0.1in
    \noindent
    \begin{lem}[Theorem 6.2.4 in \cite{B}]{\it The pull-back gerbe on $\Sigma$
    induces a canonical  element $\tilde{\theta}_{\Sigma}\in \otimes \mathcal{L}_{l_i}$, or a
    trivialization  $\theta_{\Sigma}: \otimes \mathcal{L}_{l_i} \rightarrow
    S^1$.}
    \end{lem}
    \vskip 0.1in

     Note that if we reverse the orientation of a boundary
    circle $l_i$ to obtain $\bar{l}_i$, then
    $\mathcal{L}_{\bar{l}_i}=\mathcal{L}^*_{l_i}$.

    (iii) $\theta$ has a decomposition property as follows. We decompose
    $\Sigma$ along a circle $\Sigma=\Sigma_1\cup_{S^1} \Sigma_2$
    where $\Sigma_1, \Sigma_2$ are glued along boundary circles
    $l, \bar{l}$. Let $l_l, l_{\bar{l}}$ be the corresponding
    loop. Then $l_{\bar{l}}=\bar{l}_l$ and there is a canonical
    isomorphism $\mathcal{L}_{l_l}\otimes \mathcal{L}_{l_{\bar{l}}}\cong S^1$.
    Hence $\otimes_i \mathcal{L}_{l_i}\cong \otimes_i \mathcal{L}_{l_i}\otimes
    \mathcal{L}_{l_l}\otimes \mathcal{L}_{l_{\bar{l}}}$. Under this identification,
    it is clear that
    \vskip 0.1in
    \noindent
    {\bf Gluing axiom: }{\it
    $$\theta_{\Sigma}=\theta_{\Sigma_1}\otimes
    \theta_{\Sigma_2}.$$}
    \vskip 0.1in

    $\theta$ admits another interpretation closely analogous to topological quantum field theory.
    We can view $\Sigma$ as a
    cobordism between incoming circles $l_i$ (with opposite orientation from the boundary orientation) and
    outgoing circles $l_j$.
    Then, $\theta_{\Sigma}$ can also be interpreted as
    an element of $Hom(\otimes_i \mathcal{L}_{l_i}, \otimes_j \mathcal{L}_{l_j})$. Then the gluing axiom corresponds
    to the usual gluing.

    One application of $\theta_{\Sigma}$ is to define a connection on $\mathcal{L}\rightarrow \Omega X$.
Take a
    path in the loop space
    $$F: [0,1]\times S^1\rightarrow X.$$
    Applying the above Lemma, we obtain a canonical
    isomorphism between $\mathcal{L}_{\{0\}\times S^1}$ and $\mathcal{L}_{\{1\}\times S^1}$. This can
    be viewed as the parallel transport of a connection over $\mathcal{L}$.
    Recall that a section generated by the parallel transport from
    a point is precisely a covariant constant section. From our
    construction, the restriction of a flat trivialization on
    $[0,1] \times S^1$ to each $\{t\}\times S^1$ gives a covariant
    constant section.

\vskip 0.1in
    \noindent
   \begin{lem}{\it The canonical section $s$ of $\mathcal{L}|_X$ is a covariant constant section.}
   \end{lem}
    \vskip 0.1in
    {\bf Proof: } Suppose that $F: [0,1]\times S^1\rightarrow
    X$ is a path of constant loops. Then $F$ is the
    composition of the projection to the first factor $p_1:
    [0,1]\times S^1\rightarrow [0,1]$ and a path of $X$,
    $i_F: [0,1]\rightarrow X$. We first use $i_F$ to pull back
    a gerbe with its connection to $[0,1]$. Such a pull-back is flat
    with trivial holonomy. Then, we construct a flat
    trivialization on $[0,1]$. Now, we pull it back to
    $[0,1]\times S^1$ to obtain a flat trivialization $s_F$ over
    $[0,1]\times S^1$. It is clear that the restriction of $s_F$
    to $\{t\}\times S^1$ is $s_{l_t}$ for $l_t=F(t,.)$. By the
    definition, $s_t=s_{l_t}$ is a covariant constant section along the path.

    \section{Gerbe on orbifold}
    The previous construction has been generalized to orbifolds by
    Lupercio-Uribe \cite{LU1}, \cite{LU2}. It is amazing that $\mathcal{L}\rightarrow X$
    starts to become nontrivial on an orbifold! Therefore, it is more interesting
    to study gerbes on orbifold than  on smooth manifolds!
    \subsection{Basics}

    Lupercio-Uribe's construction is carried out for an arbitrary
    groupoid. The precise definition of gerbe over an orbifold is not important
    for us. Therefore, instead of giving a long technically correct definition,
    let me motivate the definition of groupoid from orbifold. We first start
    from a smooth manifold where one can view a groupoid as a language to formalize the
    construction of an open cover. Let
    $$G_0=\bigsqcup_{\alpha} U_{\alpha}, G_1=\bigsqcup_{\alpha\beta} U_{\alpha\beta}.$$
    In the language of groupoids, $G_0$ is called the space of objects and
    $G_1$ is called the space of arrows. There are two maps
       $$s: U_{\alpha\beta}\rightarrow U_{\alpha}, t: U_{\alpha\beta}\rightarrow U_{\beta}.$$
    $s, t$ are called the source map and target map. Consider the fiber product
    $$G_2={G_1  }_{ \text{  }  t} \times_{s} G_1=\{(x, y); t(x)=s(y)\}.$$
    Using an open cover, it is not hard to see that $G_2=\bigsqcup_{\alpha\beta\gamma} U_{\alpha\beta\gamma}.$
There are also source and target maps $s, t: G_2\rightarrow G_0$
    by $s(x,y)=s(x), t(x,y)=y$. In the language of open covers, it corresponds to inclusion maps
    $$s: U_{\alpha\beta\gamma}\rightarrow U_{\alpha}, t: U_{\alpha\beta\gamma}\rightarrow
    U_{\gamma}.$$
    There is an additional multiplication map $m: G_2\rightarrow G_1$ corresponding to the
    inclusion $U_{\alpha\beta\gamma}\rightarrow U_{\alpha\gamma}$.
    To complete the definition of groupoid, we also need an
    identity $e: G_0\rightarrow G_1$ and an inverse $i:
    G_1\rightarrow G_1$. In our set-up, $e: U_{\alpha}\rightarrow
    U_{\alpha\alpha}, I: U_{\alpha\beta}\rightarrow
    U_{\beta\alpha}$ are identity maps. These structure maps
    satisfy several obvious compatibilility conditions for which we refer to
    \cite{LU1}. We often use $\G=\{G_1\stackrel{s,t}{\rightarrow}
    G_0\}$ to denote the groupoid. The process of taking a refinement of an open cover is called Morita equivalence in groupoid language.

        One can go on to construct
    $$G_{n}={G_{n-1}}_{\text{  }t} \times_{s} G_1.$$
    It corresponds to the disjoint union of $(n+1)$-fold intersections.

    With the above correspondence, we can state Lupercio-Uribe's definition of gerbes over a groupoid
    as a function $g: G_2\rightarrow S^1$ satisfying the obvious
    cocycle condition generalizing the condition on smooth
    manifolds. If two gerbes $g_1, g_2$ differ by a coboundary, we call them {\em equivalent}.
    An equivalence class of gerbes is a $\check{C}ech$ cohomology class of the sheaf $C^{\infty}(S^1)$ over the so-called classifying space
    $B\G$ of the groupoid. Furthermore, we have a long exact sequence
    $$H^2(B\G, C^{\infty}(\R))\rightarrow H^2(B\G, C^{\infty}(S^1))\stackrel{\tau}{\rightarrow} H^3(B\G, \Z).$$
    It is different from the smooth case in that the characteristic class $\tau([g])$ is an integral cohomology class of the classifying space
    $B\G$ instead of its space of orbits $|\G|$. For an orbifold, we
    can always choose its groupoid representative $\G$ with
    the property that the components of $G_0$ are contractible.
    Such a kind of groupoid is called a {\em fine} groupoid. Over
    a fine groupoid, $C^{\infty}(\R)$ is a fine sheaf. In this
    case, the equivalence class of gerbes is still classified by
    its characteristic class.

    Over a groupoid a connection on a gerbe is a pair $(A, F)$ where $A$ is a one-form on $G_1$ and
    $F$ is a two-form on $G_0$ satisfying the condition.
    $$t^*F-s^*F=dA, \ \ \ \ i \pi^*_1 A+i\pi^*_2 A +i m^*I^*A=g^{-1}dg.$$

    Now, to extend the definition of gerbe to orbifold, we just have to associate a groupoid to
    orbifolds
    which was done by Moerdijk-Pronk \cite{MP}. To this purpose, we just have to construct $G_0, G_1$ and $s,t$.

    Let $X$ be an orbifold and $(\tilde{U}_i, G_{\tilde{U}_i}), Tran(U_i, U_j)$
    be
    an orbifold atlas. We simply define $G_0=\sqcup_i \tilde{U}_i, G_1=\sqcup_{ij}Tran(U_i, U_j).$
    $s,t$ are natural projections
    $$s: Tran(U_i, U_j)\rightarrow \tilde{U}_i, \  t: Tran(U_i, U_j)\rightarrow \tilde{U}_j.$$
    We call the above groupoid {\em an orbifold groupoid}.

    An orbifold morphism corresponds to a Morita equivalence of morphisms between orbifold groupoids. An obvious and important
      fact is that
    \vskip 0.1in
    \noindent
   \begin{rem}{\it An orbifold morphism between orbifolds pulls back a gerbe with connection to a gerbe with
    connection.}
    \end{rem}
    \vskip 0.1in
    In particular, if $f: Y\rightarrow X$ is a smooth orbifold morphism from a smooth manifold $Y$ (viewed with trivial
    orbifold structure) and $[g]$ is a gerbe with connection $(A, F)$ on $X$, then $f^*[g], (f^*A, f^*F)$ is a gerbe
    with connection on a smooth manifold $Y$ even though we started from an orbifold $X$. Therefore, the previous
    construction on the holonomy line bundle $\mathcal{L}$ goes through trivially.
    However, its restriction to the inertia orbifold is no longer trivial.

   We first look at the case of discrete torsion for a global
   quotient orbifold. The inner local system has been constructed in
    \cite{R1}. We would like to show that it agrees with the holonomy
    line bundle from the gerbe induced by discrete torsion.
    Recall that discrete torsion is a two-cocycle $\alpha:
    G\times G \rightarrow S^1$. Being a cocyle means
    $$\alpha_{g,1}=\alpha_{1,g}=1, (\delta \alpha)_{g,h,k}=\alpha_{h,k}\alpha^{-1}_{gh,k}\alpha_{g,hk}\alpha^{-1}_{g,h}=1.$$
    The groupoid     presentation of the global quotient orbifold $X/G$ is a translation
    groupoid with $G_0=X, G_1=X\times G$ and
    $s(x,g)=x,t(x,g)=gx$. We will use stacky notation $[X/G]$ to denote this groupoid structure.
    One can check that $G_2=X\times G\times G$. $\alpha$
    induces a gerbe on $\G$ in the obvious way. Furthermore, we can
    choose a flat connection with $F=0, A=0$. Recall that the inertia
    orbifold is $[(\sqcup_g X^g)/G]$ and let
    $\gamma_{g,h}=\alpha_{g,h}\alpha^{-1}_{ghg^{-1},g}.$
    Recall \cite{R1} that we can define an inner local system on
    $[X/G]$ as follows. Consider the trivial bundle $\sqcup_g
    X^g \times \C_g$ where we use $\C_g$ to denote the fiber $\C$
    associated to $X^g$ and $1_g$ to denote the element
    $1$ in $\C_g$. Then, we define an action of $g: \C_h\rightarrow \C_{ghg^{-1}}$ by
    $$g(1_h)=\gamma_{g,h}1_{ghg^{-1}}.$$
    Let $\mathcal{L}_{\alpha}$ be the quotient of the trivial bundle under the
    above action.
    \vskip 0.1in
    \noindent
    \begin{thm}\label{T:discrete}{\it $$\mathcal{L}|_{\wedge [X/G]}=\mathcal{L}_{\alpha}.$$}
    \end{thm}
    \vskip 0.1in
    {\bf Proof: } We start with some algebraic preliminaries.
    If a 2-cocycle $\alpha$ can be expressed as a
    coboundary $\alpha_{g,h}=\rho_g\rho_h\rho^{-1}_{gh}$, we call
    $\rho$ a flat trivialization of $\alpha$.

    A 2-cocycle $\alpha$ corresponds to  an equivalence class of
group extensions

$$1\to S^1\to \widetilde{G}_{\alpha}\to G\to 1$$
The group $\widetilde{G}_{\alpha}$ can be given the structure of a
compact Lie group, where  $ S^1\to \widetilde{G}_{\alpha}$ is the
inclusion of a closed subgroup. The elements in the extension
group can be represented by pairs $\{ (g,a)~|~g\in G, a\in S^1\}$
with the product
$(g_1,a_1)(g_2,a_2)=(g_1g_2,\alpha_{g_1,g_2}a_1a_2)$.

    \vskip 0.1in
    \noindent
   \begin{lem}{\it There is a 1-1 correspondence between the set  of characters of
$\widetilde{G}_{\alpha}$ which restrict to scalar multiplication
on the central $S^1$ and the set of flat trivializations of
$\alpha$.}
\end{lem}
    \vskip 0.1in
    \noindent
    {\bf Proof: } If $\psi:\widetilde{G}_{\alpha}\rightarrow S^1$ is such a
character then we define an associated trivialization of $\alpha$
via $\rho (g) =\psi (g,1)$. Note that $\rho (gh)=\psi (gh,1)
=\alpha_{g,h}^{-1}\psi (gh, \alpha_{g,h})=\alpha_{g,h}^{-1} \psi
((g,1)(h,1)) =\alpha_{g,h}^{-1}\rho (g)\rho (h)$. Conversely,
given a flat trivialization $\rho :G\rightarrow S^1$, we simply
define $\psi (g,a)=a\rho (g)$.  Note that
$$\psi ((g,a)(h,b))=\psi (gh, \alpha_{g,h}ab)
=ab\rho(g)\rho(h)=a\rho (g)b\rho (h)= \psi (g,a)\psi (h,b).$$
    We have proved the lemma
    $\Box$

    Now, we come back to the proof of theorem. Any element $g\in G$
    generates an abelian subgroup $<g>$. It pulls back the 2-cocycle
    $\alpha$ and we can define the corresponding group extension
    $\widetilde{<g>}_{\alpha}$. It is well known that any 2-cocycle of a
    finite cyclic abelian group is a coboundary. Hence, we have a
    nontrivial set of flat trivializations and hence a set of
    characters of $\widetilde{<g>}_{\alpha}$. For any $h\in G$,
    $h$ sends $<g>$ to $<hgh^{-1}>$. We would like to calculate
    the action of $h$ on the set of characters (hence flat
    trivializations). Given a character $\phi$ for
    $\widetilde{<g>}_{\alpha}$ and a lifting $(h,a)$ of $h$, the
    action is defined by the formula
    $$(h,a)\phi(x,b)=\phi((h,a)(x,b)(h,a)^{-1}).$$
    $$\begin{array}{lll}
    (h,a)(x,b)(h,a)^{-1}&=&(hx, \alpha_{h,x}ab)(h^{-1},
    \alpha^{-1}_{h,h^{-1}}a^{-1})\\
    &=&(hxh^{-1},
    \alpha_{h,x}\alpha_{hx,h^{-1}}\alpha_{h,h^{-1}}^{-1}b)\\
    &=&(hxh^{-1}, \alpha_{h,x}\alpha^{-1}_{hxh^{-1},h}b)\\
    &=&(hxh^{-1},\gamma_{h,x}b)
    \end{array}$$
    Recall that $\wedge [X/G]=(\sqcup_g X^g\times \{g\})/G$. It can be
    interpreted as follows. Let $f: S^1\rightarrow [X/G]$ be a
    constant good map with Chen-Ruan characteristic $(g)$. Then,
    it factors though the constant morphism to $[X/<g>]$ which is represented by
    $(x,g)$ for $x\in X^g$.  One can
    also factor through the abelian orbifold $[X/<hgh^{-1}>]$ which
    is equivalent to the previous one. This is the action of $G$
    on $\sqcup_g X^g\times \{g\}$. Now we consider the space of
    constant morphisms to $[X/<g>]$, which is parameterized by $X^g\times \{g\}$.
    It admits a flat gerbe from $G$ through the embedding $<g>\rightarrow G$ (denoted by
    $\alpha_{<g>}$), since $<g>$ is an abelian group and such a
    gerbe is trivial. Pick a flat trivialization $\rho$ of
    $\alpha_{<g>}$. $\rho$ defines a section $s_{\rho}$ of $\mathcal{L}_{X^g\times
    \{g\}}$. We can use the same argument as in the smooth case to show
    that $s_{\rho}$ is covariant constant. Therefore, $\mathcal{L}_{X^g\times
    \{g\}}$ is trivial as a flat line bundle. However, it does not
    have a canonical flat trivialization since a different flat
    trivialization of $\alpha_{<g>}$ will define a different flat
    trivialization of $\mathcal{L}_{X^g\times  \{g\}}$. Using the
    calculation above, we conclude that the action of $G$ on
    $\mathcal{L}_{X^g\times \{g\}}$ is via the character $\gamma_{g,h}$. We
    have proved the theorem.

    \vskip 0.1in
    \noindent
    \begin{thm}[Lupercio-Uribe]{\it $\mathcal{L}|_{\wedge X}$ is flat.}
    \end{thm}
    \vskip 0.1in
    {\bf Proof: } Lupercio-Uribe proved their theorem using a
    generalization of Brylinski's relevant formula in the smooth case.
    Here, we use our previous analysis to give a direct, geometric
    proof. The  question is local. Therefore, we can assume that
    our orbifold is a global quotient $[\mathbb{R}^n/G]$ with gerbe and connection
    given by $(\alpha^0, F, A)$. It follows that $\alpha^0=\alpha \beta$ where
     $\alpha:G \times G \to S^1$ is a 2-cocyle and $\beta$ is a coboundary over $[\mathbb{R}^n/G]$ and $(\alpha^0, F, A)=(\alpha, 0, 0)+(\beta, F, A)$.
    Let $\mathcal{L}^0$,  $\mathcal{L}$ and  $\mathcal{L}'$
    represent the holonomy line bundles of $(\alpha^0, F, A)$ , $(\alpha, 0,
    0)$ and $(\beta, F, A)$ restricted to the inertia orbifoid respectively. It is clear that
     $\mathcal{L}^0 =\mathcal{L} \bigotimes\mathcal{L}'$ and ,
     fixing any path in the inertia orbifold, the parallel transport on  $\mathcal{L}^0 $ is the
     tensor product of those on $\mathcal{L}$ and  $\mathcal{L}'$.
   Since  $\beta$ is a coboundary , there is a flat trivialization
   of  $\mathcal{L}'$. The fact that  $\mathcal{L}$ is flat  implies that  $\mathcal{L}^0$ is
   flat by
   Theorem \ref{T:discrete}.

    \subsection{Holonomy on an orbifold Riemann surface
    }
       If we only consider maps from a smooth Riemann surface, the
       construction in the smooth case can be readily generalized to
       the case of orbifolds. A more interesting case is the case
       when $f: \Sigma\rightarrow X$ is a good map from an orbifold
       Riemann surface $\Sigma$. In orbifold quantum cohomology,
       we have to consider its generalization where $\Sigma$ is a
       nodal orbifold Riemann surface.

       Unfortunately, many things goes wrong and we do not have a
       straightforward generalization of string connection. One of the critical facts
       for an oriented smooth Riemann surface is $H_2(\Sigma, \Z)=\Z$ generated by its fundamental
       class $\sigma$. We use this fact to interpret the holonomy of a gerbe with a connection as
       a number in $S^1$. We
       start our discussion from the following computation of
       $H_2(B\Sigma, \Z)$ for an orbifold Riemann surface. Indeed,
       a certain subtlety arises.

    Let $\Sigma$ be an orientable orbifold Riemann surface, with
singular points $\{x_1,\cdots, x_m\}$ and corresponding
multiplicities $\{k_1,\cdots, k_m\}$. Note that the underlying
topological space $|\Sigma|$ of the orbifold $\Sigma$ is a
topological surface.  Let $B\Sigma$ be the corresponding
classifying space of the orbifold.

Given an action of the group $G$ on space $X$, let $EG$ be a free
$G$ space which is contractible. One can utilize the Borel
construction $X_G=EG\times_{G}X$. Define
$H_*^G(X,\Z)=H_*(EG\times_{G}X,\Z)$. It is well known that there
is an action of $S^1$ on a $3$-manifold $M$ such that $[M/S^1]$ is
the given orbifold $\Sigma$. It is well-known that up to a weak
homotopy equivalence $B\Sigma\cong M_G$. Hence, $H_2(B\Sigma,
\Z)=H_*^{S^1}(M, \Z)$.

There is a canonical map $\pi_X:X_G \to X/G$ and thus a
homomorphism $H_*^G(X,\Z) \to H_*(X/G,\Z)$ which is known to be an
isomorphism if the action of $G$ on $X$ is free.

 \begin{thm}{\it Let $\Sigma$ be an oriented orbifold
    Riemann surface.  \newline
 $H_2(B\Sigma,\Z)=\Z$. However, $\pi_{M*}$ does not map the generator to the fundamental class $\sigma$.
 Instead, $\pi_{M*}(e)=\pm r\sigma$ where $e\in H_2(B\Sigma, \Z)$ is the generator and $r$ is the least common multiple of the $k_i$. }
\end{thm}

    {\bf Proof: }
    All the coefficients are integers without further
    specification.

     Choose small open discs $D_i$ centered around $x_i$ such that
their closures $B_i$ are disjoint. Let $V_i=f^{-1}(B_i)$,
$V=\bigcup_iV_i$, $V_i^*=f^{-1}(D_i)$, $V^*=\bigcup_iV_i^*$ and
$U=M-V^*$. Note that the actions of $S^1$ on $U$ and $U \cap V$
are free.

Now $H_2^{S^1,\Z}(U) \cong  H_2(U/S^1,\Z) =0$ since the action is
free and $U/S^1$ is a smooth surface with boundary. $H_j(
V/S^1,\Z)=0$ for $j>0$ since  $V/S^1$ is the disjoint union of
closed discs. It is well known that the local model of  a singular
point of an orbifold Riemann surface is the quotient of a disc by
the action of a cyclic group. It follows that $H_j^{S^1,\Z}(
V)=\bigoplus_{i}H_j(B\Z_{k_i},\Z)$ and $H_2^{S^1,\Z}( V)=0$,
$H_1^{S^1,\Z}( V)=\bigoplus_i \Z_{k_i}$ since $H_2(B\Z_{k},\Z)=0$
for all positive integers.

One has the following commutative diagram by naturality of
Mayer-Vietoris sequence
    $$\begin{array}{ccccccccc}
&H_2^{S^1}(M,\Z)& \rightarrowtail &H_1^{S^1}(U \cap V,\Z)&
\rightarrow &
H_1^{S^1}(U,\Z) \oplus H_1^{S^1}( V,\Z)& \\
& \pi_{M*}\downarrow && h_1\downarrow && h_2\oplus 0\downarrow && \\
  &
H_2(M/S^1,\Z) &\rightarrowtail & H_1((U \cap V)/S^1,\Z)
&\rightarrow & H_1(U/S^1,\Z) &

    \end{array} $$

Now the homomorphisms $h_1$, $h_2$  are isomorphisms by the remark
before the statement of the Theorem. A diagram chase argument
shows that the map $p_{M*}$ is given by $a \to \pm r b$
 where $a$ is the generator of $H_2^{S^1}(M,\Z) \cong \Z$, $b$ is
  the  generator of  $H_2(M/S^1,\Z)\cong \Z$ and $r$ is the least common multiple of the $k_i$.

  It is clear that any gerbe connection $(F,A)$ on an orbifold Riemann
  surface is flat for dimension reasons. Therefore, we can define
  its holonomy $Hol(F,A)$ as a class in $H^2(B\Sigma, S^1)$ where
  $S^1$ means the constant sheaf. One can evaluate $Hol(F,A)$ on
  the generator $e\in H_2(B\Sigma, \Z)$ to obtain a number
  $Hol(F,A)(e)\in S^1$. However, the geometric evaluation is over
  the fundamental class $\sigma$. By the previous theorem,
  $p_{M*}(e)\neq \sigma$. Hence, $Hol(F,A)(\sigma)$ is undefined.

  Now, we take a different point of view of the same phenomenon.
  Consider the gerbe connection over the local 2-dimensional orbifold
  disc
  $[D/\Z_k]$. The connection as well as its holonomy is trivial. A
  choice of flat trivialization restricts to a flat trivialization
  on the boundary circle $\partial (D/\Z_k)$. Namely, we obtain an element of
  $\mathcal{L}_{\partial(D/\Z_k)}$. However, the space of flat line
  bundles
  on $[D/\Z_k]$ is non-trivial. In fact, it is parameterized by
  $\Z_k$. They induce a set of $\Z_k$-points on $\mathcal{L}_{\partial
  (D/\Z_k)}$. Now, we go back to the closed orbifold Riemann surface
  $\Sigma$ with orbifold points at $(z_1, \cdots, z_l)$ of
  multiplicity $k_1, \cdots, k_l$ and $f: \Sigma\rightarrow X$. We decompose $\Sigma$ as
  the disjoint union of orbifold discs $[D_i/\Z_{k_i}]$ and
  $V=\Sigma-\sqcup_i [D_i/\Z_{k_i}].$ Then the flat trivialization
  at $V$ specifies an element on $\mathcal{L}_{l_i}$ of each boundary
  circle $l_i$. The gluing law indicates that we should associate
  a set of $k_1\cdots k_l$-numbers in $S^1$!

    \section{Orbifold quantum cohomology twisted by a flat
    gerbe}\label{T:string}

    \subsection{String connection on an orbifold and orbifold stable maps}
     Recall the compatibility condition of an inner local system
    over each $X_{(g_1,g_2)}$;
    $\theta_{(g_1, g_2)}: e^*_1 L_{(g_1)}\otimes e^*_2 L_{(g_2)}\otimes e^*_3 L_{((g_1g_2)^{-1})}\cong 1 $.
    The purpose of this condition is as follows.
    $X_{(g_1, g_2)}=X_{(g_1, g_2, g_3)}$ with $g_3=(g_1g_2)^{-1}$ can be identified as
    the moduli space of degree zero genus zero maps with three
    marked points--$\overline{\M}_{0,3}(X, J, 0, \x)$. The evaluation maps at marked points are $e_i$.
    Let $\alpha_i \in H^*(X_{(g_i)}, L_{(g_i)})$. The trivialization
    $\theta_{(g_1, g_2)}$ maps $e_1^*\alpha_1\wedge e_2^*\alpha_2
    \wedge e_3^*\alpha_3$ to an ordinary cohomology class of
    $X_{(g_1, g_2, g_3)}$ and hence can be integrated. The latter
    property allows us to define the twisted orbifold product. To
    carry out the same construction for orbifold quantum
    cohomology, we have to construct the trivialization $\theta$
    over $\overline{\M}_{g,k}(X, J, A, \x)$ for general $g,k,A$.
    This can be accomplished by Theorem \ref{T:glue}.

    Suppose that $f: \Sigma\rightarrow X$ is an orbifold morphism. We take a
    real blow-up at the marked points to obtain a Riemann surface
    with boundary $\Sigma^{\dagger}$.
    Let  $l_{i\infty}$ be the corresponding boundary cycle. It is clear that
    each $f: \Sigma \rightarrow X$ induces a morphism $f^{\dagger}:
    \Sigma^{\dagger} \to X$. It is clear that  $f^{\dagger}_{i\infty}=f^{\dagger}(l_{i\infty})$
    is a constant loop. Moreover, we have an
    identification
    $$f_{i\infty}=e_i(f).$$
    Next, the holonomy
    $\theta_{\Sigma}=\theta_{\Sigma^{\dagger}}$ is
    interpreted as

    $$\theta_{\Sigma^{\dagger}}: \otimes_i e^*_i \mathcal{L}\rightarrow
    S^1.$$

    Next, we extend the above discussion to orbifold stable maps.
    Suppose that $f: \Sigma\rightarrow X$ is an orbifold stable
    map where $\Sigma$ is a marked orbifold nodal Riemann
    surface. For simplicity, we assume that $\Sigma=\Sigma_1\wedge
    \Sigma_2$ joining at point $p\in \Sigma_1, q\in \Sigma_2$. We
    observe that $f_{p\infty}$ is the same as $f_{q\infty}$ with
    the reversed orientation. Hence,
    $\mathcal{L}_{f_{p\infty}}=\mathcal{L}^*_{f_{q\infty}}.$
    Suppose that the marked points on $\Sigma_1$ are $x_1,\cdots,
    x_l$ and the marked points on $\Sigma_2$ are $x_{l+1},\cdots , x_k$.
    We have
    $$\theta_{\Sigma_1}: \otimes_{1\leq i\leq l} e^*_i\mathcal{L}\otimes e^*_p\mathcal{L};
    \ \ \ \ \theta_{\Sigma_2}: \otimes_{l+1\leq j\leq k}e^*_j
    \mathcal{L}\otimes e^*_q\mathcal{L}.$$
    Using the canonical isomorphism $e^*_p\mathcal{L}\otimes e^*_q\mathcal{L}\rightarrow
    S^1$, we obtain
    $$\theta_{\Sigma}: \otimes_{1\leq i\leq k}e^*_i
    \mathcal{L}\rightarrow S^1.$$
    The analogue of the gluing law is the following theorem.
    \vskip 0.1in
    \noindent
    \begin{thm}[gluing law]\label{T:glue}{\it $\theta_{\Sigma^{\dagger}}$ is continuous with respect to the
    degeneration of orbifold stable maps,i.e., it induces a
    continuous trivialization of $\otimes_i e^*_i \mathcal{L}$ over
    $\overline{\M}_{g,k}(A, X, J, \x)$.}
    \end{thm}
    \vskip 0.1in
    {\bf Proof: }
    $\theta$ is clearly continuous over each stratum of $\overline{\M}_{g,k}(A, X, J,
    \x)$. We only have to check that it is continuous with respect to
    the degeneration of orbifold stable maps. It is enough to
    discuss the case of creation of a new nodal point.
Suppose that $(f_n,(\Sigma_n,\z_n), \xi_n)$ converges to
$(f_0,(\Sigma_0,\z_0), \xi_0)$ and $p\in \Sigma_0$ is the nodal
point. It is instructive to see the degeneration of $\xi_n$ to
$\xi_0$. Recall the construction in \cite{CR2}.

Locally,
$$f_n: W_{t_n}\rightarrow (V_p-\{p\})/G_p, $$
where $W_{t}=\{xy=t; |x|, |y|<\epsilon\}$ and $(V_p, G_p,\pi_p)$
is a uniformizing system of $p\in X$. The key is to construct the
lifting $\tilde{f}_n$ mapping into $V_p$.

By Lemma 2.2.4 of \cite{CR2}, $\xi_n$ determines a characteristic
$$\theta_n: \pi_1(W_{t_n})\rightarrow G_p.$$
 Suppose that $g$ is the
image of a generator and $m$ is the order of $g$. Then $\theta_n$
determines a covering $\widetilde{W}^m_{t_n}\rightarrow W_{t_n}$.
The argument of \cite[Lemma 2.2.4]{CR2} constructs a lifting
$\tilde{f}_n: \widetilde{W}^m_{t_n}\rightarrow V_p$. Then the
convergence of $f_n$ as a good map is interpreted as convergence
of ordinary maps $\tilde{f}_n$ to $f_0:
\widetilde{W}^m_0\rightarrow V_p$, which gives a natural
compatible system $\xi_0$ at $p$. Note that $\Sigma_0$ acquires a
natural orbifold structure at the nodal point $p$, whose
uniformizing system is given by $(\widetilde{W}^m_0, \Z_m)$, where
the action of $\Z_m$ on $\widetilde{W}^m_0$ is the limit of the
action on $\widetilde{W}^m_{t_n}$.

    Let $S^1_n\subset W_{t_n}$ be the circle given by
    $|x|=\frac{|t_n|}{2\epsilon}$ with complex orientation of $x$.
    $S^1_n$ converges to a constant loop supported at $p$.  On the other hand, the same $S^1$ with
    opposite orientation can be expressed as
    $|y|=\frac{|t_n|}{2\epsilon}$ (denoted by $S^{1*}_{n}$).
    $S^{1*}_n$ converges to the constant loop support at $q$. We decompose $\Sigma_n$ along $S^1_n$ as
    $\Sigma_n=\Sigma^1_n\cup_{S^1_n} \Sigma^2_n$. Then $\Sigma^1_n$ converges to
    $\check{\Sigma}_1$ and $\Sigma^2_n$ converges to
    $\check{\Sigma}_2$. The above construction implies that $f_n|_{S^1_n}$ converges
    to $(f_0)_{p\infty}$.
    Moreover, $f_n|_{S^{1*}_n}$ as a good map converges to
    $(f_0)_{q\infty}$. . By the gluing axiom ,
    $\theta_{\Sigma_n}: \otimes_i e^*_i \mathcal{L}\rightarrow S^1$ can be
    decomposed as the product of
    $\theta_{\Sigma^1_n}: \otimes_{1\leq i\leq l} e^*_i \mathcal{L}\otimes
    \mathcal{L}_{f_n|_{S^1_n}}$ and $\theta_{\Sigma^2_n}: \otimes_{l+1\leq j\leq k} e^*_j \mathcal{L}\otimes
    \mathcal{L}_{f_n|_{S^{1*}_n}}$ using canonical trivialization
    $\mathcal{L}_{f_n|_{S^1_n}}\otimes \mathcal{L}_{f_n|_{S^{1*}_n}}$. It is clear
    that $\theta_{\Sigma_n}$ converges to $\theta_{\Sigma_0}$. $\Box$

 \subsection{Twisted orbifold Gromov-Witten invariants}
    So far, we have not yet brought in the flatness condition.
    Recall that $\otimes_i e^*_i \mathcal{L}$ in our construction is used
    as the coefficient system or flat line bundle. Therefore, we also
    need to construct a flat trivialization. This requires the
    assumption of flatness of the gerbe.
    \vskip 0.1in
    \noindent
    \begin{thm}{\it For a flat gerbe,
    $\theta_{\Sigma^{\dagger}}$ is a flat trivialization.}
    \end{thm}
    \vskip 0.1in
    {\bf Proof: } A flat bundle is completely determined by its
    holonomy around a loop. Even though
    $\overline{\M}_{g,k}(A,X,J, \x)$ is not a smooth manifold in
    general, we can still discuss the flat trivialization of a
    flat bundle. Namely, it is enough to determine if the
    trivialization is flat around each loop. Since flatness is
    a local condition, it is enough to prove that it is flat along
    a curve $f: [0,1]\rightarrow \overline{\M}_{g,k}(A,X,J, \x)$.

    Recall that there is a universal family
    $$\overline{\U}_{g,k}(A,X, J, \x)\rightarrow
    \overline{\M}_{g,k}(A,X,J,\x)$$ as an orbifold fiber bundle
    whose fiber is the domain of orbifold stable maps modulo
    automorphisms. The pullback $\pi: f^*\overline{\U}_{g,k}(A,X, J,
    \x)\rightarrow [0,1]$ is an orbifold Riemann surface bundle.
    The total space \newline $f^*\overline{\U}_{g,k}(A,X, J, \x)$ is an
    orbifold 3-manifold with boundary. Each marked point defines a
    section $s_i$. For simplicity, we first assume that there
    is no nodal point. Now, we take a real blow-up along the image
    of $s_i$ to obtain $\pi_{\dagger}: f^*\overline{\U}_{g,k}(A,X, J,
    \x)^{\dagger}\rightarrow [0,1]$. Then we replace each $s_i$ by $S^1\times [0,1]$
    (denoted by $S^1_i\times [0,1]$). It is clear that
    $\pi^{-1}_{\dagger}(t)$ is the real blow-up of $\pi^{-1}(t)$
    along the marked points,
    where $S^1_i\times \{t\}$ is precisely the infinity circle
    associated to the $i$-th marked point of $\Sigma_t=\pi^{-1}(t)$.
    A moment of thought tells that $f^*\overline{\U}_{g,k}(A,X, J,
    \x)^{\dagger}$ is an oriented 3-manifold with the boundary given
    by
    $$\partial f^*\overline{\U}_{g,k}(A,X, J,
    \x)^{\dagger}=\pi^{-1}_{\dagger}(0)\cup
    \cup_i S^1_i\times [0,1]\cup \pi^{-1}_{\dagger}(1).$$
    Furthermore, the identification happens precisely at the
    infinity circles corresponding to marked points on
    $\pi^{-1}_{\dagger}(0), \pi^{-1}_{\dagger}(1)$. Furthermore,
    there is an evaluation $e: \overline{\U}_{g,k}(A,X, J,
    \x)\rightarrow  X$ whose restriction to the $i$-th marked
    point defines $e_i$. It induces an evaluation map
    $$e_{\dagger}: f^*\overline{\U}_{g,k}(A,X, J,
    \x)^{\dagger}\rightarrow  X$$
    whose restriction to each infinity circle defines the
    corresponding evaluation map $e_{i\dagger}=f_{i\infty}$. Since the gerbe
    is flat, the holonomy around the boundary $\partial f^*\overline{\U}_{g,k}(A,X, J,
    \x)^{\dagger}$ is zero. Recall that the restriction of a flat
    trivialization of $S^1_i\times [0,1]$ to its boundary defines
    parallel transport from $\mathcal{L}_{f_{i\infty 0}}$ to
    $\mathcal{L}_{f_{i\infty 1}}$. By our definition, the restriction
    of the flat trivialization to the boundary of $\pi^{-1}_{\dagger}(0), \pi^{-1}_{\dagger}(1)$
    defines elements  $\theta_0\in \otimes \mathcal{L}_{f_{i\infty 0}}, \theta_1\in \otimes
    \mathcal{L}_{f_{i\infty 1}}$ respectively. The property that the total
    holonomy around $\partial f^*\overline{\U}_{g,k}(A,X, J,
    \x)^{\dagger}$ is zero can be interpreted as the statement
    that parallel transport maps $\theta_0$ to $\theta_1$. Then
    we prove that $\theta$ is flat.

    If we have nodal points, we do a real blow-up along the nodal
    point. It creates an additional boundary component. However, it
    is clear that the total holonomy of the additional components
    cancel each other. The above argument still applies.

    $\Box$
    \vskip 0.1in
    \noindent
    \begin{rem}{\it The proof of the above theorem depends
    critically on the flatness of the gerbe.}
    \end{rem}
    \vskip 0.1in

Now, we are ready to construct twisted orbifold GW-invariants.
Using $e_i$, we can pull back $\mathcal{L}$ (now as an orbifold
vector bundle) to define the tensor product $\otimes_i e^*_i
\mathcal{L}$. Then $\theta$ provides a flat trivialization
    $$\theta: \otimes_i e^*_i \mathcal{L}\rightarrow \C,$$
    continuous with respect to the topology of
    $\overline{\M}_{g,k}(A,X,J, \x)$. Suppose that $\x=\prod_i X_{(g_i)}$.
    It induces a homomorphism
    $$\otimes_i H^*(X_{(g_i)}, \mathcal{L}_{(g_i)})\rightarrow
    H^*(\overline{\M}_{g,k}(A,X,J,\x), \C),$$
    by
    $$(\beta_1, \cdots, \beta_k)\rightarrow \theta_*(\wedge_i
    e^*_i \alpha_i),$$
    where
    $$\theta_*: H^*(\overline{\M}_{g,k}(A,X,J,\x), \otimes
    e^*_i \mathcal{L}_{(g_i)})\rightarrow
    H^*(\overline{\M}_{g,k}(A,X,J,\x), \C)$$
    is the isomorphism induced by $\theta$.
    \vskip 0.1in
    \noindent
   \begin{defn}{\it The orbifold GW-invariants twisted by
    a flat gerbe are defined as
    $$
\Psi^{X,J}_{(g,k,A,\x, (A, F))}(\alpha^{l_1}_1, \cdots,
\alpha^{l_k}_k)= \theta_*(\prod_{i=1}^k c_1(L_i)^{l_i}e^*_i
\alpha_i )[\overline{\M}_{g,k}(X,J,A,\x)]^{vir}, $$ where $L_i$ is
the line bundle generated by the cotangent space of the $i$-th
marked point and $[\overline{\M}_{g,k}(X,J,A,\x)]^{vir}$ is the
virtual fundamental cycle constructed in \cite{CR2}.}
\end{defn}
    \vskip 0.1in
    A standard argument in Gromov-Witten theory will show that our
    twisted orbifold GW-invariants satisfy standard axioms
    \cite{CR2}. In particular, it implies that there is an
    associative quantum multiplication on $H^*_{orb}(X, \mathcal{L})$.

    \section{Computation}
    The current treatment of gerbes in the literature is usually abstract. One of the author's goals for this
    article is to be as concrete as possible. In this section, we will try to figure out how far we can go to compute
    holonomy inner local system explicitly.  We will divide this
    section into smooth versus orbifold cases.
    \subsection{Smooth case}
    When $X$ is smooth, the holonomy line bundle $\mathcal{L}|_X$ is canonically trivial. Hence $\otimes_i e^*_i \mathcal{L}|_X$
    is canonically trivial. Recall that for a stable map $f: \Sigma\rightarrow X$, $\theta_f$ is another
    trivialization $\otimes_i e^*_i\mathcal{L}|_X\cong S^1$. Therefore, the difference of two trivializations is
    a number in $S^1$. By abuse of notation, we still denote it by $\theta_f$. On the other hand, we can
    also associate the holonomy $Hol_f\in S^1$. The key theorem in the smooth case is
    \vskip 0.1in
    \noindent
   \begin{thm}{\it Suppose that $X$ is smooth and $\mathcal{L}|_X$ is trivialized by its canonical trivialization.
    Then $\theta_f=Hol_f.$}
    \end{thm}
    \vskip 0.1in
    \noindent
    {\bf Proof: }
    Let $g$ be the cocycle representing the gerbe. First, we assume that the stable map $f:
    \Sigma\rightarrow X$ has domain $\Sigma$ as an irreducible
    Riemann surface. We take a real blow-up to obtain $\Sigma^{\dagger}$.There is an obvious map $\pi:
    \Sigma^{\dagger} \rightarrow \Sigma$ by contracting $l_i$ to
    $x_i$. Let $f^{\dagger}=f\circ \pi$.
       $f^*g$ with its connection is flat on $\Sigma$. We can define its holonomy $Hol_f\in H^2(\Sigma, S^1)$.
    Since $H^2(\Sigma, S^1)=S^1$. We use $Hol_f$ to denote its $\check{C}ech$ cocycle or the number through the above
    isomorphism without any confusion. From the construction of $Hol_f$, $f^*g=Hol_f \delta
    r$. Namely, they differ by a coboundary.
    Furthermore, we are allowed to change $r$ by a constant cochain. Therefore, we can choose
    $r$ in such fashion that $r|_{\z}$ is a flat trivialization of $f^*g|_{\z}$.
    Therefore, $s_{f_i}=\pi^*r|_{x_i}$, where $f_i$ is the restriction of $f^{\dagger}$ to the boundary circle $l_i$.
    We use $f^{\dagger}$ to pull back
     the gerbe represented by the cocycle $g$. Since $H^2(\Sigma^{\dagger},
    S^1)=0$, $\pi^*Hol_f$ is trivial. Choose a trivialization $\pi^*Hol_f=\delta h$ where
    $h$ is constant. A flat trivialization of $(f^{\dagger})^*g$ is of the form
    $\delta(h \pi^*r)$. Then
    $\theta_f$ is the image of $((h\pi^*r)|_{l_1}, \cdots, (h\pi^*r)|_{l_k})$ in $\otimes_i \mathcal{L}_{f_i}$
       and $\theta_f=\prod_i h_{l_i} s_{f_i}$. We claim that $\prod_i h_{l_i}$ is
    $Hol_f$ through the canonical isomorphism $P: H^2(\Sigma, S^1)\cong S^1$. This canonical isomorphism
    is defined by the evaluation on the fundamental class of $\Sigma$. Let's consider the
    isomorphism induced by $\pi$
    $$\pi^*: H^2(\Sigma, \z, S^1)\rightarrow H^2(\Sigma^{\dagger},\partial \Sigma^{\dagger},
    S^1)$$  and the isomorphism induced by the inclusion $(\Sigma,\emptyset)\subseteq (\Sigma,\z)$ $$ H^2(\Sigma,\z, S^1)\rightarrow H^2(\Sigma, S^1).$$
     Note that $Hol_f|_{\z}=1$ and hence can be viewed as an element of $H^2(\Sigma,
    \z, S^1)$ as well as $H^2(\Sigma^{\dagger}, \partial
    \Sigma^{\dagger}, S^1)$ as its pull-back by $\pi$. Therefore, the image of $Hol_f$ under $P$ can also
    be obtained by the evaluation of $\pi^*Hol_f$ on the relative fundamental class of
    $(\Sigma^{\dagger}, \partial \Sigma^{\dagger})$ which  is precisely $\prod_i h_{l_i}$ since
    $<\pi^*Hol_f,[\Sigma^{\dagger},\partial \Sigma^{\dagger}]>=<\delta h,[\Sigma^{\dagger},\partial\Sigma^{\dagger}]>
    =<h,[\partial \Sigma^{\dagger}]>= \prod_i h_{l_i}$ .

    Now, we consider the general case in which $\Sigma$ may have more than one component.
    Then, we apply previous argument for each component. By the gluing axiom, $\theta_f$ is multiplicative.
    $Hol_f$ is obviously multiplicative. Hence, the theorem is true for multi-components $\Sigma$
    as well.
 $\Box$

    Suppose that we have a flat connection $(A, F)$. Then we have the global
    holonomy $Hol\in H^2(X, S^1)$. It is clear that $Hol_f=Hol(f_*[\Sigma])$. Then we prove that
    \vskip 0.1in
    \noindent
   \begin{cor}\label{T:smooth}{\it Suppose that $(g,F,A)$ is a flat gerbe. Under the canonical trivialization of $\mathcal{L}|_X$,
     the twisted GW-invariant
    $$\Psi^{X,J}_{(g,k,A,\x, (A, F))}(\alpha^{l_1}_1, \cdots,
\alpha^{l_k}_k)=Hol(A) \Psi^{X,J}_{(g,k,A,\x)}(\alpha^{l_1}_1,
\cdots, \alpha^{l_k}_k).$$
    }
    \end{cor}
    \vskip 0.1in

    \subsection{Discrete torsion}

    Suppose that we have a global quotient orbifold $[X/G]$ and a
    discrete torsion $\alpha: G\times G\rightarrow S^1$. Recall
    that we can express the inertia orbifold as a global quotient
    $\wedge [X/G]=[(\sqcup_g X^g\times \{g\})/G]$. Furthermore,
    the holonomy line bundle
    $$\mathcal{L}_{\wedge [X/G]}=(\sqcup_g X^g\times
    \{g\})\times_{\gamma}\C.$$
    We would like to express the moduli space of orbifold stable
    maps as a $G$-global quotient as well. This has already been
    done in \cite{JKK}. Let's briefly review their construction.
    As we mentioned previously, by pulling back via the $G$-bundle
    $X\rightarrow X/G$, an orbifold stable map is equivalent to a
    $G$-orbifold cover $E\rightarrow \Sigma$ and a $G$-map $\phi:
    E\rightarrow X$. A $G$-stable map has additional data
    $\tilde{z}_i$-a lifting of the marked point $z_i\in \Sigma$.
    Consider the isotropy subgroup $G_{\tilde{z}_i}\subset G$. For dimensional
    reasons
    $G_{\tilde{z}_i}\cong \Z_k$ for some $k$. $\tilde{z}_i$ uniquely determines an element $g_i$ as the
    generator of $G_{\tilde{z}_i}$. $g_i$ has a geometric
    interpretation as the monodromy of a small loop around $z_i$. It
    is independent of the lifting of small loop because a
    different lifting will conjugate $g_i$ by an element of
    $G_{\tilde{z}_i}$, which is abelian. Now, the evaluation map
    $e_i$ is lifted to
    $$\tilde{e}_i(E\rightarrow \Sigma, \phi, \tilde{z}_1, \cdots
    \tilde{z}_k)=(\phi(\tilde{z}_i), g_i).$$
    $G$ acts on $G$-stable maps by its action on $\tilde{z}_i$ and
    $\tilde{e}_i$ is $G$-equivariant. Let
    $\overline{\M}^G_{g,k}(A,X,J)$ be the moduli space of
    $G$-stable maps. We can apply virtual fundamental cycle
    techniques for $\overline{\M}^G_{g,k}$ to obtain a $G$-orbifold
    GW-invariant. The orbifold GW-invariant is the invariant part
    of the $G$-orbifold GW-invariant.

    We can use $\tilde{e}_i$ to pull back $\mathcal{L}_{\wedge [X/G]}$ to
    form a flat line bundle
    $$\overline{\M}^G_{g,k}(A,X,J)\times_{\gamma_{g_1}\cdots\gamma_{g_k}}
   ( \otimes_i \C_{g_i}).$$
    By our construction the holonomy of $\Sigma$ should give a
    flat trivialization of this flat line bundle and we would like
    to write it down explicitly.

    Without loss of generality, we assume that $\Sigma$ is
    irreducible. We take a real blow-up $\Sigma^{\dagger}$ at all
    the marked points. There results a real blow-up $E^{\dagger}$ of $E$ at all the
    preimage points of $z_i$. $G$ acts on $E^{\dagger}$ freely and
    $\Sigma^{\dagger}=E^{\dagger}/G$. Using the translation groupoid
    representation $[E^{\dagger}/G]$ of $\Sigma^{\dagger}$, we can pull back $\alpha$ to
    a flat gerbe on $\Sigma^{\dagger}$. The general theory tells
    us that this flat gerbe is trivial. Furthermore, a flat
    trivialization restricts to a flat trivialization on each
    boundary circle. Let's study the induced gerbe on
    $\Sigma^{\dagger}$ from the point of view of the Chen-Ruan
    characteristic. Fix a base point $x_0$ in the interior of
    $\Sigma^{\dagger}$ and choose a lifting $\tilde{x}_0\in
    E^{\dagger}$. It defines the CR-characteristic $\rho:
    \pi_1(\Sigma^{\dagger}, x_0)\rightarrow G$. In fact,
    $E^{\dagger}=E_{univ}\times_{\rho} G$ where $E_{univ}$ is the universal cover. We can use $\rho$ to
    pull back $\alpha$ to a 2-cocycle $\tilde{\alpha}$ of
    $\pi_1(\Sigma^{\dagger}, x_0)$. The induced gerbe on
    $\Sigma^{\dagger}$ is given by $\tilde{\alpha}$.  Since
    $E_{univ}$ is contractible,
    $$H^2(\pi_1(\Sigma^{\dagger}, x_0), S^1)\cong H^2(\Sigma^{\dagger},
    S^1)\cong 0.$$
    Therefore, $\tilde{\alpha}$ is a coboundary. Choose  a
    cobundary $\tilde{\alpha}=\delta h$.

    Let's go back to the local monodromy $g_i$ at $\tilde{z}_i$.
    Indeed, it is the monodromy along the boundary circle associated with
    $z_i$. We would like to embed $<g_i>$ in
    $\pi_1(\Sigma^{\dagger}, x_0)$ as a subgroup. This can be done
    as follows. Choose a path $d_i$ from $\tilde{x}_0$ to the
    boundary circle associated with $\tilde{z}_i$ with the end point $x_i$,
    then go around the boundary circle to $g_i x_i$ and go back to
    $g_i \tilde{x}_0$ along $g_i d^{-1}_i$. Its projection
    $l_i$ on $\Sigma^{\dagger}$ is a loop at $x_0$ whose
    lifting defines $\rho(l_i)=g_i$. Then, we map $g_i$ to
    $\gamma_i$. It can be shown that a different path $d'_i$ conjugates $l_i$ by
    an element of the image of $\pi_1(E^{\dagger}, \tilde{x}_0)$. The image of
    $\pi_1(E^{\dagger}, \tilde{x}_0)$ is precisely the kernel of $\rho$.
    It is clear that $\tilde{\alpha}_{<g_i>}=\alpha_{<g_i>}=\delta
    h_{<g_i>}$. Therefore, we can use $h_{<g_i>}$ to trivialize
    $\mathcal{L}_{X^{g_i}\times \{g_i\}}$. Hence, the trivialization given by holonomy on $\Sigma^{\dagger}$
    corresponds to the pull-back of the trivial bundle $\mathcal{L}_{X^{g_i}\times \{g_i\}}$.
    We still have to show that $\gamma_{g_1}\cdots \gamma_{g_k}=1$
    in order to descend to a trivial bundle over
    $\overline{\M}_{g,k}(A, [X/G], J)$. This follows directly from
    the CR-characteristic $\rho$. Recall that $\pi_1(\Sigma^{\dagger},
    x_0)$ has generators $l_1, \cdots, l_k, \mu_1,
    \cdots \mu_{2g}$ with relation $\prod_j [\mu_{2j-1},
    2j] l_1 \cdots l_k=1$. Therefore,
    $$1=\gamma_{\prod_j [\mu_{2j-1},
    \mu_{2j}]l_1\cdots l_k}=\prod_j[\gamma_{\mu_{2j-1}},
    \gamma_{\mu_{2j}}]\gamma_{l_1}\cdots\gamma_{l_k}=\gamma_{l_1}\cdots
    \gamma_{l_k}.$$
    Then we can apply the previous construction to the $G$-orbifold Gromov-Witten invariants constructed in
    \cite{JKK}. Thus, we proved
    \vskip 0.1in
    \noindent
  \begin{thm}{\it The twisted orbifold Gromov-Witten invariant
    of a discrete torsion is the $G$-invariant part of the
    $G$-orbifold Gromov-Witten invariant, under the action twisted
    by the discrete torsion.}
    \end{thm}

{\bf Address:}

\begin{description}
    \item[Jianzhong Pan] Institute of Mathematics, Academica Sinica

   \item[Yongbin Ruan] Department of Mathematics, Hong Kong University of Science and Technology and University of
      Wisconsin-Madison

     \item[Xiaoqin Yin]     Department of Mathematics, Hong Kong University of Science and Technology
\end{description}
    \end{document}